\documentclass[a4paper,twoside,12pt]{article}

\usepackage{verbatim}
\usepackage{mathrsfs}
\usepackage[french]{babel}
\usepackage[T1]{fontenc}
\usepackage{amsmath}
\usepackage{amssymb}
\usepackage{amsthm}
\usepackage{amscd}
\usepackage{relsize,exscale}
\usepackage{hyperref}
\usepackage{cleveref}

\usepackage{tabularx}
\usepackage{latexsym}
\usepackage{amsfonts}
\usepackage{longtable}
\usepackage{epsfig}
\usepackage{graphicx}
\usepackage{color}
\usepackage{float}
\usepackage{multicol}
\usepackage{hyperref}
\addto\captionsfrench{}
\addto\captionsfrench{}
\addto\captionsfrench{}
\newtheorem{thm}{Theorem}[section]

\newtheorem{prop}[thm]{Proposition}

\newtheorem{defn}[thm]{Definition}
\newtheorem{example}[thm]{Example}

\newtheorem{prob}[thm]{Problem}
\newtheorem{rem}[thm]{Remark}

\newenvironment{dem}{{\noindent {\bf Proof. }}}{\hfill {\rule{2mm}{2mm}} }

\newcommand{\symm}{\mathfrak{S}}

\newcounter{chapitre}
\newcommand\rdots{\mathinner{\mkern1mu\raise0pt\vbox{\kern7pt\hbox{.}}
     \mkern2mu\raise4pt\hbox{.}\mkern2mu\raise8pt\hbox{.}\mkern1mu}}

\def\N{\mathcal{N}}

\definecolor{gris25}{gray}{0.75}
\definecolor{gris20}{gray}{0.80}

\def\keywords#1{{ \it{Keywords:}} {#1}} %

\date{}
\begin{document}

%\maketitle
%\mainmatter  % start of an individual contribution

% first the title is needed
\title{Width-$k$ Eulerian polynomials of type $A$ and $B$:\\  the $\gamma$-positivity}

\maketitle

\vspace{.1cm}

\begin{center}
\author{Marwa Ben Abdelmaksoud$^{1}$\;  and Adel Hamdi$^{2}$}

\end{center}
\maketitle
\begin{center}

\centerline{\footnotesize{$^{1,2}$Faculty of Science of Gabes, Department of Mathematics,}} \centerline{\footnotesize{ Cit\'e Erriadh 6072,  Zrig, Gabes, Tunisia}}
%\centerline{\footnotesize{and}}
%\centerline{\footnotesize $^{}$Institut Camille Jordan, Universit\'e de Lyon, Universit\'e Lyon 1,} \centerline{\footnotesize
%UMR 5208 du CNRS, 43, boulevard du 11 novembre 1918}
%\centerline{\footnotesize 69622, Villeurbanne Cedex, France}
\centerline{\footnotesize \texttt{$^{1}$abdelmaksoud.marwa@hotmail.fr \quad $^{2}$aadel\_hamdi@yahoo.fr.}}
\end{center}

\begin{abstract}In this paper, we introduce some new generalizations of classical descent and inversion statistics on signed permutations that arise
from the work of Sack and \'Ulfarsson \cite{13}, and called
\emph{$k$-width descents} and \emph{$k$-width inversions} of type
$A$ (\cite{2}). Using the aforementioned new statistics, we derive
new generalizations of Eulerian polynomials of type $A$,
    $B$ and $D$. We establish also the $\gamma$-positivity of the Eulerian \emph{"width-k"} polynomials. Referring to Petersen's paper \cite{10}, we give a
     combinatorial interpretation of finite sequences associated with these new polynomials using quasi-symmetric functions and a partition $P$.
\end{abstract}

 % Write down at least 3 Keywords

%\tableofcontents
\section{ Introduction}
The main purpose of this paper is to extend some fundamental aspects of the theory of Eulerian polynomials on Coxeter groups and their unimodality, symmetry and $Gamma$-positivity. Many polynomials with combinatorial meanings have been shown to be unimodal (see \cite{4} or \cite{5} for example). Let $\mathcal{A}=\{a_i\}_{i=0}^{d}$ be a finite sequence of nonnegative numbers. Recall that a polynomial $g(x)=\sum_{i=0}^{d} a_i x^{i}$ of degree $d$ is said to be positive and unimodal, if the coefficients are increasing and then decreasing, i.e., there is a certain index $0\leq j \leq d $ such that $$0\leq a_0\leq a_1\leq\cdots \leq a_{j-1}\leq a_j \geq a_{j+1} \geq\cdots \geq a_d \geq 0.$$

We will say that $g(x)$ is $palindromic$ ( or symmetric) with center of symmetry at $\lfloor d/2\rfloor$, if $a_i=a_{d-i}$ for all $0\leq i\leq d.$\\
The polynomial $g(x)$ is said to be $Gamma$-positive (or $\gamma$-positive) if $$g(x)=\sum_{i=0}^{\lfloor d/2\rfloor} \gamma_i x^{i} (1+x)^{d-2 i},$$
where $d$ is the degree of $g(x)$ and nonnegative reals $\gamma_0,\gamma_1,\ldots,\gamma_{\lfloor d/2\rfloor}$.
So, both palindromic and unimodal are two necessary conditions for the $Gamma$-positivity of $g$($x$). One of the most important polynomials in combinatorics is the nth Eulerian polynomials, for the statistic $"des^A"$ on $\symm_n$, defined as
$$A_n(x)= \displaystyle\sum_{\sigma \in \symm_n} x^{des^A(\sigma)}.$$
Given a set of combinatorial objects $\tau$, a combinatorial statistics is an integer for to every element of the set. In other words, it is a function $st :\tau\rightarrow \mathbb{N}$.\\
For a statistic $st$ on symmetric group $\symm_n$, one may form the generating function: $$F^{st}_n(x)=\displaystyle \sum_{\pi \in \symm_n} x^{st(\pi)}.$$

MacMahon \cite{1} considered four different statistics for a permutation $\sigma$ in the group of all permutations $\symm_n$ (it is also called type-$A$ permutation) of the set $[n]:=\{1, \ldots, n\}$; the number of descents $(des^A(\sigma))$, the number of excedances $(exc^A(\sigma))$, the number of inversions $(inv^A(\sigma))$ and the major index $(maj^A(\sigma))$. Given a permutation $\sigma$ in $\symm_n$, we say that the pair $(i,j)\in [n]^{2}$ is an inversion of $\sigma$ if $i<j$ and $\sigma(i)>\sigma(j)$, that $i\in \{1, 2, \ldots, n-1\}$  is a descent if $\sigma(i)>\sigma(i+1)$, and that $i \in \{1, \ldots, n\}$ is an excedance if $\sigma(i)>i$. The major index is the sum of all the descents. These four statistics have many generalization.

The symmetric group $\symm_n$ is generated by the set $\tau:=\{\tau_i;~ 1 \leq i \leq n-1\}$, where $\tau_i:=(1,\ldots, i-1, i+1, i, i+2, \ldots, n) \in \symm_n$, for $i \in [n-1]$. Therefore, the length function on $ \symm_n$ is defined to be: for any $\sigma \in \symm_n$,
$$\ell^{A}(\sigma):=\textnormal{min}\{r\geq 0;~ \sigma=\tau_{x_1}\tau_{x_2}\cdots \tau_{x_r};~ \tau_{x_i} \in \tau,~ x_i \in [n-1] \}.$$

Then the length function can be written in terms of inversion number on $\symm_n$  (see \cite{14}, Proposition 3.1, for $\sigma \in \symm_n$ ) as follows
  $$\ell^{A}(\sigma)=inv^A(\sigma).$$

 The paper is organized as follows. We start with some definitions which generalized the \emph{width-k descents} and \emph{width-k inversions} statistics on classical permutations studied by Davis \cite{2} into signed permutations. In section $2$, we will prove Proposition \ref{lb}, in which we improve the combinatorial formulas of these last statistics in signed permutations and give some examples. In section $3$ and $4$, we will show Theorem \ref{gammaA}, Theorem \ref{thmb} and Theorem \ref{gammaB}, in which we define the \emph{width-k Eulerian} polynomials of type $A$ and $B$. So, we give some recurrence relations concerning the coefficients of these polynomials. Then, we will study the $\gamma$-positivity by specifying the combinatorial values of $\gamma$. Finally, in section $5$, which is the same as section $3$ and $4$, we will define the \emph{width-k Eulerian} polynomials of type $D$ and we define two sets $WD_{n,k,p}$ and $W\bar{D}_{n,k,p}$ in order to find the recurrence relations for the coefficients of this polynomial. We will prove Theorem \ref{gammaD} by studying the necessary condition for this polynomial to be $\gamma$-positive.
\section{Width-k descents and width-k inversions on signed permutations}\label{section2}
Recently, Sack and \'Ulfarsson \cite{13} introduced new
generalizations of classical descents and inversions statistics for
any permutation in $\symm_n$. They are called \emph{width-$k$
descents} and \emph{width-$k$ inversions} (see \cite{2}). For each
$1\leq k < n$, the sets of such statistics are defined as follows
\begin{eqnarray*}
Des^{A}_k(\sigma) &:=&\{i \in [n-k]; \sigma(i)>\sigma(i+k)\},\\
Inv^{A}_k(\sigma)&:=&\{(i,j)\in[n]^{2};\sigma(i)>\sigma(j)~~and~j-i=mk,~m>0
\}.
\end{eqnarray*}
Their cardinalities are denoted, respectively, by
$des^{A}_k(\sigma)$ and $inv^{A}_k(\sigma)$.\\

In this paper, we study same analogues of these statistics on signed
permutations. A signed permutation is a bijection of $[-n,n]:=\{-n,
\cdots, -1, 1, \cdots, n\}$ onto itself that satisfies
$\pi(-i)=-\pi(i)$ for all $i \in [n]$. We denote by $B_n$ the set of
signed permutations of length $n$, which are known as
hyperoctahedral groups. Let $D_n\subset B_n$ be the subset
consisting of the signed permutations with even negative entries
number. Adin, Brenti and Roichman \cite{3} defined a permutation
statistics known as signed descent number (or type-$B$ descent
number) and $flag$ descent number. A signed descent of $\pi=(\pi(1),
\pi(2), \ldots, \pi(n))\in B_n$ is an integer $0\leq i\leq n-1$
satisfying $\pi(i)>\pi(i+1)$, where $\pi(0)=0$. The signed descent
of $\pi \in D_n$ (type-$D$ descent) given by Chow \cite{9} is a
type-$B$ descent restricted to $D_n$. The signed descent number of
$\pi \in B_n$ is denoted by $des^B$($\pi$) and defined as follows
$$des^B(\pi):=|\{0\leq i\leq n-1; ~\pi(i)>\pi(i+1)\}|.$$
The $flag$ descent statistics of a signed permutation $\pi$, denoted
by $fdes^B(\pi)$, counts a descent at position $0$ once and all
other descents twice. In other words,
$$fdes^B(\pi) :=des^B(\pi)+des^A(\pi).$$
The \emph{nth Eulerian} polynomials of signed permutations and the
\emph{nth flag descents} polynomials are defined, respectively, by
$$\mathfrak{B}_n(x)=\displaystyle \sum_{\pi \in B_n}
x^{des^B(\pi)},$$
 $$F_n(x)=\displaystyle \sum_{\pi \in B_n}
x^{fdes^B(\pi)}.$$
\begin{defn}
Let $\pi$ be a permutation in $B_n$. We define:
\begin{enumerate}
\item the inversion number in $\pi$ by
 $$inv^A(\pi):=|\{(i,j); 1 \leq i < j \leq n ~\textnormal{and}~\pi(i)>\pi(j)\}|,$$
\item the negative integers number in $\pi$ by
$$neg(\pi):=|\{i \in [n];\pi(i)<0\}|,$$
\item the negative pairs sum number on $\pi$ by $$nsp(\pi):=|\{(i,j)\in [n]^{2}; i<j ~\textnormal{and}~ \pi(i)+\pi(j)<0 \}|,$$
\item the type-$B$ inversion's number of $\pi$ by $$inv^B(\pi):=inv^A(\pi)+neg(\pi)+nsp(\pi).$$
\end{enumerate}
\end{defn}
As the Coxeter group, $B_n$ is generated by the set $\tau^B:=\{\tau^B_1,~\tau^B_2, \ldots, \tau^B_{n-1},~\tau^B_0\}$, where $\tau^B_0:=(-1, 2, \cdots, n)$ and $\tau^B_i:=(1,2, \cdots, i-1, i+1, i,i+2, \cdots, n)$, for $i \in [n-1]$ (for more details see again \cite{14}). We define the length function $\ell^{B}(\cdot)$, on $ B_n$, by
$$\ell^{B}(\pi):=\textnormal{min}\{r\geq 0;~ \pi=\tau^B_{x_1}\tau^B_{x_2}\cdots \tau^B_{x_r};~ \tau^B_{x_i} \in \tau^B,~ 0 \leq x_i \leq n-1\},\quad \forall\pi \in B_n.$$

\noindent Recently, Brenti (\cite{14}) has proved that the function $\ell^{B}(\cdot)$ can be interpreted combinatorially as
 $$\ell^{B}(\pi)=inv^B(\pi).$$

\noindent We extend the same definitions of the width-$k$ statistics on classical permutations, introduced in \cite{2}, to signed permutations as follows:
\begin{defn}
Let $1\leq k \leq n$ and $\pi=(\pi(1),\,\pi(2),\cdots,\,\pi(n))$ be a permutation in $B_n$ with $\pi(0) =0$. We define:
\begin{enumerate}
\item $Inv^{A}_{k}(\pi):=\{(i,j); 1 \leq i < j \leq n;\pi(i)>\pi(j)~\textnormal{and}~j-i=mk,m>0\}$ as the set of width-$k$ type $A$ inversion,
\item $Inv^B_k(\pi):=\{(i,j); 0\leq i<j \leq n;~\pi(i)>\pi(j)~\textnormal{and}~j-i=mk,m>0\}\cup\{(-i,j); 1\leq i<j \leq n;~\pi(-i)>\pi(j)~\textnormal{and}~j-i=mk,m>0\}$ as the set of width-$k$ inversion,
\item $Des^{B}_{k}(\pi):=\{0\leq i \leq n-k;~\pi(i)>\pi(i+k)\}$ as the set of width-$k$ descent,
\item $des^{B}_{k}(\pi):=|Des^{B}_{k}(\pi)|$ as the width-$k$ descent number,
\item $fdes^{B}_k :=des^{A}_{k}(\pi)+ des^{B}_{k}(\pi)$ as the width-$k$ flag descent number,
\item $ndes^{B}_k(\pi):=|\{1 \leq i \leq n-k;\pi(-i)>\pi(i+k)\}|$ as the width-$k$ negative descent number,
\item $inv^{A}_{k}(\pi):=|Inv^{A}_{k}(\pi)|$ as the width-$k$ type $A$ inversion number,
\item $neg_k(\pi):=|\{1 \leq i \leq \lfloor \frac{n}{k}\rfloor; ~\pi(ik)<0\}|$ as the width-$k$ negative number,
\item $nsp_k(\pi):=|\{(i,j)\in [n]^{2};\pi(i)+\pi(j)<0, j-i=mk,m>0\}|$ as the width-$k$ negative sum pairs number,
\item $inv^B_k(\pi):= inv^A_k(\pi)+neg_k(\pi) + nsp_k(\pi)$ as the width-$k$ inversion number.
\end{enumerate}
\end{defn}
\noindent For the case $k=1$, the set $Inv^B_1(\pi)$ can be also denoted by $Inv^B(\pi)$ and we will prove that $inv^B_k(\pi) = |Inv^B_k(\pi)|$, for all $1\leq k \leq n$.\\

\noindent Let $des^{D}_{k}(\pi)$ be the set of all \emph{width-$k$} descents of $\pi \in D_n$ given by
$$des^{D}_{k}(\pi):=|\{0\leq i \leq n-k;\pi(i)>\pi(i+k)\}|,~~(\pi(0)=0).$$
Taking $k=1$, we recover Borowiec's and  Mlotkowski's definition of width-$1$ type $D$ descents (see \cite{12} for $\pi(0) =0$). A similar definition was established by Brenti (\cite{14}), in which $\pi$(0):=$-\pi$(2).\\

\noindent Let $K\subseteq [n]$ ($\neq\emptyset$)  be the widths set under consideration. It holds that:
$$Des^{B}_{K}(\pi):=\displaystyle\bigcup_{k\in K}Des^{B}_{k}(\pi)~\textnormal{and}~des^{B}_{K}(\pi)=|Des^{B}_{K}(\pi)|,$$
$$Inv^{A}_{K}(\pi):=\displaystyle \bigcup_{k\in K}Inv^{A}_{k}(\pi)~\textnormal{and}~inv^{A}_{K}(\pi)=|Inv^{A}_{K}(\pi)|,$$
$$Inv^{B}_K(\pi):=\displaystyle\bigcup_{k\in K}Inv^B_k(\pi)~\textnormal{and}~ inv^B_K(\pi) = |Inv^{B}_K(\pi)|.$$
\begin{example}\label{ff}
For $\pi=(4,-1,-3,-6,5,-7,2)\in B_7$, we obtain:
\begin{enumerate}
\item $Des^{B}_{\{2,3\}}(\pi)=\{0,1,2,3,4,5\}$,
\item $Inv^{A}_{\{2,3\}}(\pi)=\{(1,3),(1,4),(1,7),(2,4),(2,6),(3,6),(4,6),(5,7)\}$,
\item $Inv^{B}_{\{2,3\}}(\pi)=\{(0,2),(0,4),(0,6),(1,3), (1,7), (2,4),(2,6),(3,7),(5,7), (-2,4),$\\ $(-2,6),(-3,7),(-5,7),(0,3), (1,4), (3,6), (-1,4), (-3,6), (-4,7)\}$.
\end{enumerate}
\end{example}
In the following two propositions, we generalize the width-$k$ inversion number of type $A$ definition established in \cite{2}.
\begin{prop}\label{new prop1} For any  $k \in [n]$ and any $\pi \in B_n$, we have
\begin{equation}\label{eq6}
|Inv^{B}_{k}(\pi)| = inv^{B}_{k}(\pi)=\displaystyle\sum_{m\geq 1} (des^{B}_{mk}(\pi)+ndes^{B}_{mk}(\pi)).
\end{equation}
\end{prop}
\begin{dem}
For some $m>0$, the elements of $Inv^{B}_{k}(\pi)$ are pairs of the form $(i,i+mk)$; $0\leq i \leq n-1$ or $(-i,i+mk)$; $1 \leq i \leq n-1$. So, an element exists in $Inv^{B}_{k}(\pi)$ if and only if there is a width-$mk$ descent of $\pi$ at $i$ or a width-$mk$ negative descent of $\pi$ at $(-i)$. Thus, $inv^{B}_{k}(\pi)$ just counts the number of descents of length $mk$ for every possible $m$. Hence, the identity (\ref{eq6}) holds true.

\end{dem}
\begin{prop}\label{new prop2}
For any $K\subseteq [n]$ and any $\pi \in B_n$, we have
\begin{equation}\label{eq8}
inv^{B}_{K}(\pi)=\displaystyle\sum_{\emptyset\subsetneq K'\subseteq K}(-1)^{|K'|+1}~inv^{B}_{lcm(K')}(\pi),
\end{equation}
where $inv^{B}_{lcm(K')}(\pi)=0$ for $lcm(K')\geq n+1.$
\end{prop}
\begin{dem} If $K = \{k\}$ then the identity (\ref{eq6}) allows to conclude.
For a general subset $K$, applying the same technics developed in \cite{2}-Proposition $2$, for the classical permutations, we can conclude.

\end{dem}
\begin{example}
Let $\pi$ as in Example \ref{ff}. Then, we have
\begin{align*}
inv^{B}_{2}(\pi)&=\displaystyle\sum_{m\geq 1} (des^{B}_{2m}(\pi)+ndes^{B}_{2m})\\
&=(des^{B}_{2}(\pi)+ndes^{B}_2(\pi))+(des^{B}_{4}(\pi)+ndes^{B}_4(\pi))+(des^{B}_{6}(\pi)+ndes^{B}_6(\pi))\\
&=(5+2)+(2+2)+(2+0)=13.
\end{align*}
Moreover, we have $inv^{B}_{\{2,3\}}(\pi)= 19$, where $inv^{B}_2=13$ and $inv^{B}_3=8$. Remark that $(0,6)$ and $(1,7)$ have both the width $2$ and the width $3$, so it must also have the width-$lcm(2,3)$. Thus, $inv^{B}_{6}(\pi)=\{(0,6),(1,7)\}.$ Hence, it holds that $inv^{B}_{\{2,3\}}(\pi)=inv^{B}_{2}(\pi)+inv^{B}_{3}(\pi)-inv^{B}_{6}(\pi).$
\end{example}
We aim, now, to generalize the search function described in \cite{2}, on the set of signed permutations. It helps to show the interaction between width-$k$ statistics by changing its normalization map.\\

\noindent Let $n,\, k$ be positive integers satisfying $n=dk+r$, for some $(d,r) \in \mathbb{N}^2$ with $0\leq r <k$. For any $\pi$ in $B_n$, we may associate the set of disjoint substrings $\gamma_{n,k}(\pi)=\{\gamma_{n,k}^{1}(\pi),\gamma_{n,k}^{2}(\pi), \ldots, \gamma_{n,k}^{k}(\pi)\}$, where
\[\gamma_{n,k}^{i}(\pi)=\begin{cases} (\pi(i), \pi(i+k), \pi(i+2k), \ldots, \pi(i+dk)),& \text{if $i\leq r;$}\\
(\pi(i), \pi(i+k), \pi(i+2k), \ldots, \pi(i+(d-1)k)),& \text{if $r<i\leq k.$} \end{cases}\]
Let $\psi$ be the map defined as follows
\begin{eqnarray}\label{r1}
\psi:B_n & \longrightarrow& B^{r}_{d+1}\times B^{k-r}_{d}\nonumber\\
 \pi &\longmapsto&(std \gamma_{n,k}^{1}(\pi),std \gamma_{n,k}^{2}(\pi), \ldots, std \gamma_{n,k}^{k}(\pi)).
 \end{eqnarray}
We denote by $std$ the standardization map. For all $1 \leq i \leq k$, the permutation $std \gamma_{n,k}^{i}(\pi)$ is obtained by replacing
 the smallest integer in absolute value of $\gamma_{n,k}^{i}(\pi)$ by $1$, the second smallest integer in absolute value by $2$, etc...
  Then, for each element of $\gamma_{n,k}^{i}(\pi)$, we add $"-"$ at each $\pi(i+jk)<0$, where $0 \leq j \leq
  d$. This yields that each $std\gamma_{n,k}^{i}(\pi)$ is a signed
permutation of $B_d$ or $B_{d+1}$.
\begin{example}
Let $\pi$ as in Example \ref{ff} and assume that $k=3$. We have
\begin{align*}
\gamma_{7,3}(\pi)&=(std \gamma_{7,3}^{1}(\pi),std \gamma_{7,3}^{2}(\pi),std \gamma_{7,3}^{3}(\pi))\\
&=(std(4,-6,2),std(-1,5),std(-3,-7))\\
&=((2,-3,1),(-1,2),(-1,-2)).
\end{align*}
\end{example}

Let $\ell^{B}_k(.)$ be the \emph{width-k}-analogue definition of the length function statistics on signed permutations, given by
$$\ell^{B}_k(\pi):= \displaystyle\sum_{i=1}^k \ell^{B}(std \gamma_{n,k}^{i}(\pi)).$$
In the following result, we give the explicit combinatorial description of $\ell^{B}_k(.)$.
\begin{prop} Let $\pi \in B_n$. Then, we have
$$ \ell^{B}_k(\pi)=inv^{B}_k(\pi).$$
\end{prop}
\begin{dem} For any $\pi\in B_n$, it is straightforward to see that
\begin{eqnarray*}
\left|Inv_k^B(\pi)\right|&=&\left|\displaystyle \bigcup_{i=1}^k Inv^B(std \gamma_{n,k}^{i}(\pi))\right|\\
&=&\displaystyle \sum_{i=1}^k \left|Inv^B(std \gamma_{n,k}^{i}(\pi))\right|\\
&=&\displaystyle \sum_{i=1}^k \ell^B(std \gamma_{n,k}^{i}(\pi))= \ell^{B}_k(\pi).
\end{eqnarray*}
\end{dem}

For $n,\,k\in \mathbb{N},\;n=dk+r$ with $0\leq r<k,\;d>0$, we denote by $M_{n,k}$ the multinomial coefficient defined by
$$M_{n,k}=\dbinom{n}{(d+1)^{r}, d^{k-r}},$$
where $i^m$ indicates $i$ repeated $m$ times.
Let $[ n]_{x}$ be the $x$-analogue of the integer $n\geq 1$ given by
$$[ n]_{x}:=\frac{1-x^{n}}{1-x}=1+x+x^{2}+\cdots +x^{n-1}.$$
Obviously, we obtain the $x$-analogue factorial
$$[n]_{x}!:=[1]_{x}\cdots ~[n-1]_{x}[n]_{x},~~[0]!:=1~~ \textnormal{and}~~ [2n]_x!!:=\prod_{i=1}^{n} [2i]_x,~~[0]!!:=1.$$
We come now to one of the main results in this framework.
\begin{thm}\label{lb} For $n \geq 1$, we have
\begin{equation}\label{nsp}
\displaystyle \sum_{\pi \in B_n} x^{inv^A(\pi)+nsp(\pi)} t^{neg(\pi)}=[n]_x! \displaystyle \prod_{i=0}^{n-1}(1+t x^{i}).
\end{equation}
\end{thm}
\begin{dem}
We use an induction argument to show (\ref{nsp}). Notice that the result is obvious for $n=1$.
Assume that the result holds true for the step $n-1$ and we will show it for $n$ ($n\geq 1$), that is
\begin{equation}\label{p(n-1)}
\displaystyle \sum_{\pi \in B_{n-1}} x^{inv^A(\pi)+nsp(\pi)} t^{neg(\pi)}=[n-1]_x! \displaystyle \prod_{i=0}^{n-2}(1+t x^{i}).
 \end{equation}
For that, it suffices to take into account the number of
\emph{inversions}, \emph{nsp} and \emph{negatives} for any integer
$\pi(n)$ in $\pi=(\pi(1) , \pi(2), \ldots, \pi(n-1),\pi(n))\in B_n$.
So, if $\pi(n)=l$, for $1 \leq l \leq n $ and for all $j \in
[-n,n]$, such that $\pi(j)=k$ where $ l< k \leq n$, we have exactly
one inversion or one nsp. Then $\pi(n)$ makes $(n-l)$ choices of
\emph{inversions} and \emph{nsp}. This implies that the identity
\eqref{p(n-1)} is multiplied by $x^{n-l}$, for all $\pi(n)=l$. If
$\pi(n)=-l$, $1 \leq l \leq n$, then for all $-l < k < l$ there
exists $j \in [-n,n]$, such that $\ pi(j)=k$. In this case, we have
exactly $(2 l -2)$ \emph{inversions} and \emph{nsp}. And for all $l
< k \leq n$ if $k>l$, we have exactly $(n-l)$ \emph{inversions} or
\emph{nsp}. Thus, $\pi(n)$ makes $(n-2+l)$ choices of
\emph{inversions} and \emph{nsp}. This gives the identity
\eqref{p(n-1)} multiplied by $t x^{n-2+l}$. Finally, the identity
\eqref{p(n-1)} is multiplied by $[n]_x(1+tx^{n-1})$. This completes
the proof.

\end{dem}

\noindent As an immediate consequence of Theorem \ref{lb}, for
$t=1$, the identity \eqref{nsp} becomes
\begin{equation}\label{nsp2}
\displaystyle \frac{1}{2}\sum_{\pi \in B_n} x^{inv^A(\pi)+nsp(\pi)}=
[n]_x [2(n-1)]_x!!\,.
\end{equation}
and we recover the similar identity as in the OEIS [\cite{28},
A162206].
\section{Width-k Eulerian polynomials of type $A$ }
In this section, we will introduce Davis's \emph{width-k Eulerian
polynomials} of type $A$ and we will prove its
\emph{$\gamma$-positivity}. Thus, the concept of
\emph{$\gamma$-positivity} on classical Eulerian polynomials,
    one of the most important polynomials in combinatorics, appeared first in the works of Foata and Sch\"{u}tzenberger \cite{6} and then of Foata
     and Strehl $(\cite{7},\cite{8})$. \emph{$\gamma$-positivity} is an elementary property that polynomials with symmetric coefficients can have,
      which directly implies their unimodality. Working on the statistics of \emph{Eulerian descent} is in a way a
generalization of the study of eulerian numbers which counts the
number of permutations with the same descent number. For a
permutation $\sigma \in \symm_n$, an index $i \in [n]$ is a double
descent of $\sigma$ if $\sigma(i-1)> \sigma(i)> \sigma(i+1)$, where
$\sigma(0)=\sigma(n+1)=\infty$. We also have a \emph{left peak}
(resp. \emph{pic}) of $\sigma \in \symm_n$ is any index $i \in
[n-1]$ (resp. $2 \leq i \leq n-1$) such that
$\sigma(i-1)<\sigma(i)>\sigma(i+1)$, where $\sigma(0):=0$.\\

\noindent In the following, we define such statistics.
\begin{defn} Let $\sigma=(\sigma(1), \sigma(2), \ldots, \sigma(n))$ be a permutation in $\symm_n$ and $1\leq k < n$.
\begin{enumerate}
  \item The set of double width-k descents in $\sigma$ is given by
  $$Ddes^{A}_{k}(\sigma):=\{i \in [n]; \sigma(i-k)> \sigma(i)> \sigma(i+k)\}$$
   where $\sigma(j)=\infty$, for all $j>n$ or $j\leq 0$. The number of double width-$k$ descents in $\sigma$ is
   $$ddes^{A}_{k}(\sigma):=|Ddes^{A}_{k}(\sigma)|.$$
  \item The set of width-$k$ peaks $($called also interior width-$k$ peaks$)$ in $\sigma$ is given by
   $$Peak_k(\sigma):=\{ k+1\leq i \leq n-k ; \sigma(i-k)<\sigma(i)>\sigma(i+k)\}.$$
    The number of width-$k$ peaks is
    $$peak_k(\sigma):=|Peak_k(\sigma)|.$$
  \item The set of width-$k$ left peaks in $\sigma$ is given by
   $$Lpeak_k(\sigma):=\{k \leq i \leq n-k ; \sigma(i-k)<\sigma(i)>\sigma(i+k)\}$$ where $\sigma(0)
   =0.$
        The number of width-$k$ left peaks in $\sigma$ is
         $$lpeak_k(\sigma):=|Lpeak_k(\sigma)|.$$
         \end{enumerate}
\end{defn}
\noindent The \emph{width-$k$ Eulerian} polynomials of type $A$,
denoted by $\mathfrak{WA}_{n,k}(x)$, are with the following form
\begin{equation}
\mathfrak{WA}_{n,k}(x):=  F^{des^{A}_k}_n(x)= \displaystyle \sum_{\sigma \in \symm_n} x^{des^{A}_k(\sigma)}.
\end{equation}
We denote by $\mathrm{WA}_{n,k,p}$ the set
$$\mathrm{WA}_{n,k,p}:=\{\sigma \in \symm_n; des^{A}_k(\sigma)=p\},$$ and $a(n,k,p)$ its
cardinal. Taking $k = 1$, we recover the classical Eulerian
polynomials, and its $nth$ $\gamma$-positivity,
$A_n(x)=\mathfrak{WA}_{n,1}(x)$, given by Foata and
Sch\"{u}tzenberger [\cite{6}, Theorem 5.6], with form
%$A_{n,1}(x)$ was first shown combinatorially by Foata and Sch\"{u}tzenberger $(\cite{6},~ Theorem 5.6)$.
\begin{equation}
\mathfrak{WA}_{n,1}(x)=\displaystyle\sum_{p=0}^{\lfloor
\frac{n-1}{2}\rfloor} |\Gamma_{  n,p}|~ x^{p} (1+x)^{n-1-2 p},
\end{equation}
where $\Gamma_{ n,p}$ is the set of permutations $\sigma \in S_{n}$
with $des^{A}_{1}(\sigma)=p$ and $ddes^{A}_{1}(\sigma)=0$. In the
same context, an other interpretation was given by Shapiro et al.
\cite{11} by introducing the notion of slides, and showing that the
$\gamma$-expansion coefficient $\gamma_{n,p}$ counts the number of
$n$-permutations with $p$ descents and $p$ slides. Later, Chow
generalized in \cite{16} the notion of slides to types $B$ and $D$.

\noindent The following proposition gives an identity that was
originally established by Sack and \'Ulfarsson \cite{13}, but with a
different notation. Later, Davis \cite{2} gives a different proof.
\begin{prop} For all $n\geq1 ~and ~1\leq k\leq n-1,$ we have
\begin{equation}\label{desa}
\mathfrak{WA}_{n,k}(x)= F_{n}^{des^{A}_{k}}(x)= M_{n,k} A_{d+1}^{r}(x) A_{d}^{k-r}(x).
\end{equation}
\end{prop}
\noindent Let $\alpha$($n,k,p$) be the coefficients of the
polynomials $\mathfrak{WA}_{n,k}(x)$ such that $a(n,k,p) = M_{n,k}
\alpha$($n,k,p$). For a clearer observation, we give, in the table
below, the first coefficients of $\alpha$($n,k,p$),
  for $1 \leq n \leq 6,~1\leq k \leq n-1$ and $0 \leq p \leq
  n-k$.
  \newpage
\begin{table}[!h]
\centering
\begin{tabular}{|c|c|c c c c c c|}
\hline
\scriptsize n & \scriptsize k &
\multicolumn{6}{c|}{\centering{\scriptsize p}}\\ [-.1cm]
\cline{3-8}
& & \scriptsize0 &\scriptsize1 &\scriptsize2 & \scriptsize3 &\scriptsize4 &\scriptsize 5   \\ [-.1cm]
\hline
\scriptsize1&\scriptsize1&\scriptsize1& &  & & &\\ [-.1cm]
\hline
\scriptsize2&\scriptsize1&\scriptsize1&\scriptsize1& & & &\\ [-.1cm]
\hline
\scriptsize3&\scriptsize1&\scriptsize1&\scriptsize4&\scriptsize1& & &\\[-.1cm]
 &\scriptsize2&\scriptsize1&\scriptsize1& & & &\\[-.1cm]
\hline
\scriptsize4&\scriptsize1&\scriptsize1&\scriptsize11&\scriptsize11&\scriptsize1& &\\[-.1cm]
 &\scriptsize2&\scriptsize1&\scriptsize2&\scriptsize1& & &\\[-.1cm]
 &\scriptsize3&\scriptsize1&\scriptsize1 & & & &\\[-.1cm]
\hline
\scriptsize5&\scriptsize1&\scriptsize1&\scriptsize26&\scriptsize66&\scriptsize26&\scriptsize1&\\[-.1cm]
&\scriptsize2&\scriptsize1&\scriptsize5&\scriptsize5&\scriptsize1& &\\[-.1cm]
&\scriptsize3&\scriptsize1&\scriptsize2&\scriptsize1& & &\\[-.1cm]
&\scriptsize4&\scriptsize1&\scriptsize1& & & &\\[-.1cm]
\hline
\scriptsize6&\scriptsize1&\scriptsize1&\scriptsize57&\scriptsize302&\scriptsize302&\scriptsize57&\scriptsize1\\[-.1cm]
&\scriptsize2&\scriptsize1&\scriptsize8&\scriptsize18&\scriptsize8&\scriptsize1&\\[-.1cm]
&\scriptsize3&\scriptsize1&\scriptsize3&\scriptsize3&\scriptsize1& &\\[-.1cm]
&\scriptsize4&\scriptsize1&\scriptsize2&\scriptsize1& & &\\[-.1cm]
&\scriptsize5&\scriptsize1&\scriptsize1& & & &\\
\hline
\end{tabular}
\caption{The first values of $\alpha(n,k,p)$.}
\end{table}

The polynomial $\mathfrak{WA}_{n,k}(x)$ (as the product of $k$ unimodal, symmetric and $\gamma$-positive polynomials) is unimodal, symmetric with
nonnegative coefficients, $\gamma$-positive with
 center of symmetry $\lfloor \frac{n-k}{2}\rfloor$ and $deg(\mathfrak{WA}_{n,k}(x))=n-k$.
 % Because, it is the product of $k$ unimodal, symmetric and $\gamma$-positive polynomials.\\
For example,
$$\mathfrak{WA}_{6,2}(x)=M_{6,2}
(1+8x+18x^{2}+8x^{3}+x^{4})$$ is $\gamma$-positive since
$$\mathfrak{WA}_{6,2}(x)=20 x^{0}(1+x)^{4}+80x(1+x)^{2}+80x^{2}(1+x)^{0}.$$
The coefficient $\alpha(n,k,p)$ is a new  integers sequence on the
Coxeter group of type $A$. It is therefore natural to pose the
following problem on the recurrence relation of this sequence which
has been confirmed by the fact that any $n \geq 4$, $k$ is the
smallest positive integer such that $n=2k+r $ with $0\leq r <k$ and
$0 \leq p \leq n-k$, we have
\begin{eqnarray*}
\alpha(n,k,p)&=&\alpha(n-2,k-1,p-1)+\alpha(n-2,k-1,p),\\
\hbox{with}\quad\alpha(n,k,0)&=&\alpha(n,k,n-k)=1\quad\textnormal{and}\quad
\alpha(n,k,-1)=0. \end{eqnarray*}
\begin{prob}
Is it possible to find the recurrence relation of $\alpha(n,k,p)$, for all $1 \leq k \leq n$?
\end{prob}
We have the following result.
\begin{thm}\label{gammaA} For all $n\geq1 ~and ~1\leq k\leq n-1,$ we have
$$\mathfrak{WA}_{n,k}(x)=\displaystyle \sum_{\sigma \in \symm_n} x^{des^{A}_k(\sigma)}=\displaystyle\sum_{p=0}^{\lfloor \frac{n-k}{2}\rfloor} |\Gamma_{ n,k,p}|~ x^{p} (1+x)^{n-k-2 p},$$
where $\Gamma_{ n,k,p}$ is the set of permutations $\sigma \in S_{n}$ with $des^{A}_{k}(\sigma)=p$ and $ddes^{A}_{k}(\sigma)=0$.
\end{thm}

\begin{dem} Using the standardization map $\varphi$ defined on $\symm_n$ (see \cite{2}) and assuming that, for all $1 \leq k \leq n-1$, $\varphi(\sigma)=(\sigma_1, \sigma_2, \ldots, \sigma_k)$ such that $\sigma_i = std \gamma_{n,k}^{i}(\sigma)$ for all $i$, we have $$des^{A}_{k}(\sigma)=\displaystyle \sum_{i=1}^{k} des^A(\sigma_i).$$
So, it is clear that each \emph{width-k} descent and \emph{width-k} double descent in $\sigma$ are usual descent and double descent in some unique $\sigma_i$ and vice versa. Since $\mathfrak{WA}_{n,k}(x)$ is a $\gamma$-positive with center of symmetry $\lfloor \frac{n-k}{2}\rfloor$ and $deg(\mathfrak{WA}_{n,k}(x))=n-k$, the desired claim follows from that.

\end{dem}

\noindent For $1 \leq n \leq 6,~1\leq k \leq n-1$ and $0 \leq p \leq \lfloor \frac{n-k}{2} \rfloor$, we give in the following tabular some values of $\gamma^{A}_{ n,k,p}.$
\begin{table}[!h]
\centering
\begin{tabular}{|p{.1cm}|p{.1cm}|p{.3cm} p{.3cm} p{.3cm}|}
\hline
\scriptsize n & \scriptsize k &
\multicolumn{3}{c|}{\centering{\scriptsize p}}\\[-.1cm]
\cline{3-5}
& &\scriptsize 0&\scriptsize1&\scriptsize2\\[-.1cm]
\hline
\scriptsize1&\scriptsize1&\scriptsize1& &\\[-.1cm]
\hline
\scriptsize2&\scriptsize1&\scriptsize1& &\\[-.1cm]
\hline
\scriptsize3&\scriptsize1&\scriptsize1&\scriptsize2 &\\[-.1cm]
 &\scriptsize2&\scriptsize3& &\\[-.1cm]
\hline
\scriptsize4&\scriptsize1&\scriptsize1&\scriptsize8&\\[-.1cm]
 &\scriptsize2&\scriptsize6&\scriptsize0&\\[-.1cm]
 &\scriptsize3&\scriptsize12& &\\[-.1cm]
 \hline
\scriptsize5&\scriptsize1&\scriptsize1&\scriptsize22&\scriptsize16\\[-.1cm]
 &\scriptsize2&\scriptsize10&\scriptsize20&\\[-.1cm]
 &\scriptsize3&\scriptsize30&\scriptsize0&\\[-.1cm]
 &\scriptsize4&\scriptsize60& &\\[-.1cm]
\hline
\scriptsize6&\scriptsize1&\scriptsize1&\scriptsize52&\scriptsize136\\[-.1cm]
 &\scriptsize2&\scriptsize20&\scriptsize80&\scriptsize80\\[-.1cm]
 &\scriptsize3&\scriptsize90&\scriptsize0&\\[-.1cm]
 &\scriptsize4&\scriptsize180&\scriptsize0&\\[-.1cm]
 &\scriptsize5&\scriptsize360& & \\
\hline
\end{tabular}
\caption{The first values of $\gamma^{A}_{  n,k,p}$.}
\end{table}

\section{Width-k Eulerian polynomials of type $B$}
In this section, we give a new generalization of type $B$ Eulerian polynomials and its $\gamma$-positivities. We introduce $\mathfrak{WB}_{n,k}(x)$ the \emph{width-k Eulerian} polynomials of type $B$ as follows
 $$\mathfrak{WB}_{n,k}(x)= F^{des^{B}_k}_n(x):=\displaystyle \sum_{\pi \in B_n} x^{des^{B}_k(\pi)}.$$
Let $\mathrm{WB}_{n,k,p}$ be the set
$$\mathrm{WB}_{n,k,p}:=\{\pi \in B_n; des^{B}_k(\pi)=p\}.$$
Its cardinal will be denoted by $b(n,k,p)$. It is clear that, for $k = 1$, we obtain the classical Eulerian polynomials of type $B$. The associated $nth$ $\gamma$-positive is given by the following result.
\begin{thm} $(\cite{10},\textnormal{Proposition 4.15})$
 For all $n\geq1$, we have
  $$\mathfrak{WB}_{n,1}(x)=\displaystyle\sum_{p=0}^{\lfloor n/2\rfloor} \gamma^{B}_{ {n,p}}~ x^{p} (1+x)^{n-2p},$$
   where $\gamma^{B}_{  {n,p}}$ is the number of permutations $\sigma \in \symm_n$ with $p$ left peaks, multiplied by $4^{p}$.
\end{thm}

For any $\pi \in B_n$, we define the statistics descent of type $A$ over the set of signed permutations as follows
$$des^A(\pi):=|\{i \in [n-1]; \pi(i)>\pi(i+1)\}|.$$
The following identities holds true.
\begin{thm}\label{thmb} For any $n \geq k \geq 1$ and $d\geq 0$ such that $n = dk+r$, $0\leq r<k$, we have
\begin{eqnarray}\label{desb}
F_{n}^{des^{B}_{k}}(x)&=& 2^{n-d} M_{n,k} B_{d}(x) A_{d}^{k-r-1}(x) A_{d+1}^{r}(x),\\
F_{n}^{fdes^{B}_{k}}(x)&=&2^{n-d} M_{n,k}~F_d(x)~A_{d}^{k-r-1}(x^{2})~ A_{d+1}^{r}(x^{2}),\\
F^{inv^{B}_k}_n(x)&=&2^{k-1} M_{n,k}~ [d]^{k-r-1}_x ~[d+1]^{r}_x~ [2d]^{r+1}_x!! ~[2(d-1)]^{k-r-1}_x!! .
\end{eqnarray}
\end{thm}

\begin{dem} Let $\psi$ be the map given by (\ref{r1}).
%$$\psi:B_n\rightarrow B^{r}_{d+1}\times B^{k-r}_{d}$$
%$$\psi(\pi)=(std \gamma_{n,k}^{1}(\pi),std \gamma_{n,k}^{2}(\pi), \ldots, std \gamma_{n,k}^{k}(\pi)).$$
We fix $k\in[n]$ and $\psi(\pi)=(\pi_1, \pi_2, \ldots, \pi_k)\in B^{r}_{d+1}\times B^{k-r}_{d}$ such that, $std \gamma_{n,k}^{i}(\pi)=\pi_i$ for all $i$. There are $M_{n,k}$ choices to partition $[n]$ into subsequences $\gamma_{n,k}^{i}(\pi)$.
For  $\pi \in B_n$, we define the map $\epsilon$ as follows
\[\epsilon(\pi)=\begin{cases} 1, ~~\text{if $\pi(1)<0;$}\\ 0, ~~\text{otherwise.}\end{cases}\]
We note that $$des^{B}_k(\pi)=des^B(\pi_k)+\displaystyle\sum_{i=1}^{k-1} (des^B(\pi_i)-\epsilon(\pi_i)).$$
It follows from that
 $$\displaystyle\sum_{i=1}^{k-1} (des^B(\pi_i)-\epsilon(\pi_i))=\displaystyle\sum_{i=1}^{k-1} des^A(\pi_i),\quad \pi \in B_n.$$
Hence,
\begin{align*}
F^{des^{B}_k}_n(x) &=\displaystyle\sum_{\pi \in B_n} x^{des^{B}_k}(\pi)\\
&=M_{n,k} \sum_{\pi_k \in B_d,~(\pi_1, \ldots, \pi_{k-1})\in B^{k-r-1}_d\times B^{r}_{d+1}} x^{des^B(\pi_k)} x^{des^A(\pi_1)}\cdots x^{des^A(\pi_{k-1})}\\
&=M_{n,k} ~B_d(x) ~(2^{d})^{k-r-1}A^{k-r-1}_d(x) ~(2^{d+1})^{r} ~A^{r}_{d+1}(x)\\
&=2^{n-d}~ M_{n,k} ~B_d(x)~ A^{k-r-1}_d(x) ~A^{r}_{d+1}(x)
\end{align*}
which proves the first identity.

\noindent Now, by using the definition of width-$k$ descent, we get
\begin{align*}
fdes^{B}_k(\pi)&=des^{A}_k(\pi)+des^{B}_k(\pi)\\
&=\sum_{i=1}^{k}des^A(\pi_i)+des^B(\pi_k)+\sum_{i=1}^{k-1}des^A(\pi_i)\\
&=\sum_{i=1}^{k-1}2 des^A(\pi_i)+des^A(\pi_k)+des^B(\pi_k)\\
&=fdes^B(\pi_k)+\sum_{i=1}^{k-1}2 des^A(\pi_i).
\end{align*}
On the other hand, we have $$\displaystyle \sum_{\pi \in B_n} x^{2des^A(\pi)}= 2^{n} A_n(x^{2}).$$
Hence,
\begin{align*}
F^{fdes^{B}_k}_n(x)&=\sum_{\pi \in B_n} x^{fdes^{B}_k(\pi)}\\
&=M_{n,k} \sum_{\pi_k \in B_d~,(\pi_1, \ldots, \pi_{k-1})\in B^{k-r-1}_d\times B^{r}_{d+1}} x^{fdes^B(\pi_k)} x^{2 des^A(\pi_1)}\cdots x^{2 des^A(\pi_{k-1})}\\
&=2^{n-d}~ M_{n,k}~ F_d(x)~  A^{k-r-1}_d(x^{2}) ~A^{r}_{d+1}(x^{2}).
\end{align*}
which proves the second identity. Now, observe that the width-$k$ inversion number in $\pi \in B_n$ is the sum
$$inv^{B}_k(\pi)=inv^{B}(\pi_k)+\displaystyle \sum_{i=1}^{k-1}(inv^A(\pi_i)+nsp(\pi_i)).$$
The generating function of the inversion numbers can be presented as following (see, for
instance \cite{15}, Section $3.15$)
\[\displaystyle\sum_{\pi \in B_n} x^{inv^B(\pi)} =[2n]_x!!\,.\]
Hence, we obtain
\begin{align*}
F^{inv^{B}_k}_n(x) &=\displaystyle\sum_{\pi \in B_n} x^{inv^{B}}_k(\pi)\\
&=M_{n,k}\sum_{\pi_k \in B_d,~(\pi_1, \ldots,\pi_{k-1})\in B^{k-r-1}_d\times B^{r}_{d+1}} x^{inv^{B}(\pi_k)} x^{(inv^A(\pi_1)+nsp(\pi_1))}\cdots x^{(inv^A(\pi_{k-1})+nsp(\pi_{k-1}))}\\
&=2^{k-1} M_{n,k} [2d]^{r+1}_x!! [2(d-1)]^{k-r-1}_x!! [d]^{k-r-1}_x [d+1]^{r}_x
\end{align*}
which improves the last identity.

\end{dem}

By using the identity (\ref{desb}), we deduce that
$$\mathfrak{WB}_{n,k}(x)=2^{n-d} M_{n,k} B_{d}(x) A_{d}^{k-r-1}(x)
A_{d+1}^{r}(x).$$
 Let $\beta$($n,k,p$) be the coefficient of the polynomial $\mathfrak{WB}_{n,k}(x)$ such that $$b(n,k,p)=2^{n-d}M_{n,k}\beta(n,k,p).$$
In the table below, we give some coefficients of $\beta$($n,k,p$),
for $1 \leq n \leq 6, 1\leq k \leq n$ and $0 \leq p \leq n-k+1$.
\begin{table}[!h]
\centering
\begin{tabular}{|c|c|c c c cc c c|}
\hline
\scriptsize n & \scriptsize k &
\multicolumn{7}{c|}{\centering{\scriptsize p}}\\
\cline{3-9}
& & \scriptsize0 &\scriptsize1 &\scriptsize2 & \scriptsize3 &\scriptsize4 & \scriptsize5 &\scriptsize6\\
\hline
\scriptsize1&\scriptsize1&\scriptsize1&\scriptsize1 &  & & & & \\
\hline
\scriptsize2&\scriptsize1&\scriptsize1&\scriptsize6 &\scriptsize1  & & & & \\[-.2cm]
 &\scriptsize2&\scriptsize1&\scriptsize1 &  & & & & \\
\hline
\scriptsize3&\scriptsize1&\scriptsize1&\scriptsize23 &\scriptsize 23 &\scriptsize 1& & & \\[-.2cm]
 &\scriptsize2&\scriptsize1&\scriptsize2 & \scriptsize1 & & & & \\[-.2cm]
 &\scriptsize3&\scriptsize1&\scriptsize1 &  & & & & \\
\hline
\scriptsize4&\scriptsize1&\scriptsize1&\scriptsize76 & \scriptsize230 &\scriptsize 76&\scriptsize1 & & \\[-.2cm]
 &\scriptsize2&\scriptsize1&\scriptsize7 & \scriptsize7 &\scriptsize1 & & & \\[-.2cm]
 &\scriptsize3&\scriptsize1&\scriptsize2 & \scriptsize1 & & & & \\[-.2cm]
 &\scriptsize4&\scriptsize1&\scriptsize1 &  & & & & \\
\hline
\scriptsize5&\scriptsize1&\scriptsize1&\scriptsize237 &\scriptsize1682  &\scriptsize1682 & \scriptsize237& \scriptsize1& \\[-.2cm]
 &\scriptsize2&\scriptsize1&\scriptsize10 & \scriptsize26 & \scriptsize10&\scriptsize 1& & \\[-.2cm]
 &\scriptsize3&\scriptsize1&\scriptsize3 & \scriptsize3 &\scriptsize 1& & & \\[-.2cm]
 &\scriptsize4&\scriptsize1&\scriptsize2 & \scriptsize 1& & & & \\[-.2cm]
 &\scriptsize5&\scriptsize1&\scriptsize1 &  & & & & \\
\hline
\scriptsize6&\scriptsize1&\scriptsize1&\scriptsize722 &\scriptsize10543  &\scriptsize23548 &\scriptsize10543 &\scriptsize722 & \scriptsize1\\[-.2cm]
 &\scriptsize2&\scriptsize1&\scriptsize27 &\scriptsize116  &\scriptsize 116& \scriptsize27& \scriptsize1& \\[-.2cm]
 &\scriptsize3&\scriptsize1&\scriptsize8 &\scriptsize 14 &\scriptsize8 &\scriptsize 1& & \\[-.2cm]
 &\scriptsize4&\scriptsize1&\scriptsize3& \scriptsize3 &\scriptsize1 & & & \\[-.2cm]
 &\scriptsize5&\scriptsize1&\scriptsize2 &\scriptsize 1 & & & & \\[-.2cm]
 &\scriptsize6&\scriptsize1&\scriptsize1 &  & & & & \\
\hline
\end{tabular}
\caption{The first values of $\beta(n,k,p)$.}
\end{table}

Since the coefficient $\beta(n,k,p)$ is a new integer sequence on the Coxeter group of type $B$,
it is natural to pose the following problem on the recurrence relation of this sequence.
 It has been confirmed that any $n \geq 3$, $k$ is the smallest positive integer such that $n+1=2k+r$ with $0\leq r <k$ and $0 \leq p \leq n- k+1$,
  we have the following recurrence relation
\begin{eqnarray*}
\beta(n,k,p)&=&\beta(n-2,k-1,p-1)+\beta(n-2,k-1,p),\\
\hbox{with}\quad
\beta(n,k,0)&=&\beta(n,k,n-k+1)=1,\quad\textnormal{and}\quad\beta(n,k,-1)=0.
\end{eqnarray*}
\begin{prob}
Is it possible to find the recurrence relation of $\beta(n,k,p)$, for all $1 \leq k \leq n$?
\end{prob}

The polynomial $\mathfrak{WB}_{n,k}(x)$ is unimodal and symmetric with nonnegative coefficients. It is $\gamma$-positive (as the product of $k$ $\gamma$-positive polynomials) with center of symmetry $\lfloor \frac{n-k+1}{2}\rfloor $ and $deg(\mathfrak{WB}_{n,k}(x))=n-k+1$. For example, for $n=6,\,k=2$, we have $$\mathfrak{WB}_{6,2}(x)=2^{3} M_{6,2} (1+27x+116x^{2}+116x^{3}+27x^{4}+x^{5})$$ is $\gamma$-positive since
$$\mathfrak{WB}_{6,2}(x)=160 x^{0}(1+x)^{5}+3520x(1+x)^{3}+6400x^{2}(1+x).$$
So, we have the following theorem.
\begin{thm}\label{gammaB} For any $1\leq k\leq n,$ the following identity holds true
$$\mathfrak{WB}_{n,k}(x)=\displaystyle \sum_{\pi \in B_n} x^{des^{B}_k(\pi)}=\displaystyle\sum_{p=0}^{\lfloor \frac{n-k+1}{2}\rfloor} 2^{2p+k-1} |\Gamma^{\mathcal{(\ell)}}_{  n,k,p}|~ x^{p} (1+x)^{n-k+1-2 p},$$
where $\Gamma^{\mathcal{(\ell)}}_{ n,k,p}$ is the set of permutations $\sigma \in S_{n}$ with $p$ width-$k$ left peaks.
\end{thm}
\noindent To prove this theorem, we need to generalize some results
on $P$-partitions and enriched $P$-partitions due to the work of
Petersen \cite{10}.\\

In the partition theory, it is well known that there are two
$P$-partition definitions. The first one is due to Stanley \cite{20}
which defined them as the order reversing maps. The second one is
introduced by Gessel \cite{21} and defined them as order preserving
maps. In the current paper, we will adopt the second definition.\\

\noindent Let $X$ be a subset of the positive integers and $\mathfrak{L}(P)$ be the set of all permutations of $[n]$ which extend a poset
 $P$ with partial order $<_p$ to a total order. When $X$ has a finite cardinality $p$, the number of $P$-partitions must also be finite.
  In this case, let define the order polynomial, denoted by $\Omega(P;p)$, to be the number of $P$-partitions $f : [n] \to X$.
   In our study of $P$-partitions, it is enough to consider $P$ as a permutation and $\Omega_A(P,p)$ is the type $A$ order
   polynomial. Recall that Chow \cite{23} defined a type $B$ poset $P$ whose elements are 0, $\pm 1$, $\ldots$,
   $\pm n$ such that if $i<_p j$ then $-j<_p -i$.
 So, the type $B$ $P$-partitions differ from ordinary $P$-partition only in the property $f(-i) =-f(i)$.

\noindent For the group of signed permutations, the only difference
with the symmetric group of type $A$ is that if $\pi(1)<0$,
 then $0$ is a descent of $\pi$. Let $\Omega_B(\pi,p)$ be the order polynomial for any signed permutation.
 For fixed $n$, the polynomials of order turn out to be: $\Omega_A(i,x)=\Omega_B(i,x)=\binom{x+n-i}{n}$,
 for any permutation of type $A$ or type $B$ with $i-1$ descents. Similarly to corollary $2.4$ of \cite{10}, we will give the relation
between the polynomial of order width-$k$ of a poset $P$ and the sum
of the polynomials of order width-k of its linear extensions in the
following definition.
\begin{defn}
We can say that each permutation $\pi \in \symm_n$ as a poset with the total order $\pi(s) < \pi(s+k)$, for all $s \in Des_k^A(\pi)$. Then, the width-k order polynomial of a poset P having $p \in \N$ descents, denoted by $\Omega_A(P,k,p)$, is the sum of the width-k order polynomials of its linear extensions:$$\Omega_A(P,k,p)=\displaystyle\sum_{\pi \in\mathfrak{L}(P) }\Omega_A(\pi,k,p),$$ where, for any permutation $\pi\in \symm_n$, the width-k order polynomial is
\begin{equation}
\Omega_A(\pi,k,p):=| \{f:[n-k+1] \to [p]/
\end{equation}
\begin{gather*}
1 \leq f(\pi(1))\leq f(\pi(1+k))\leq f(\pi(1+2 k))\leq \cdots f(\pi(1+\tau_d k))< \\ f(\pi(2))\leq f(\pi(2+k))\leq f(\pi(2+2 k))\leq \cdots f(\pi(2+\tau_d k))<\cdots < \\f(\pi(k))\leq f(\pi(k+k))\leq f(\pi(k+2k))\leq \cdots \leq f(\pi(k+\tau_d k)) \leq p,\\ and ~f(\pi(s))< f(\pi(s+k)),~ if~ s\in des^{A}_k(\pi)\}|,
\end{gather*}
with \[\tau_d=\begin{cases} d, &\text{if$ ~i\leq r $;}\\
d-1, &\text{if ~$r< i \leq k$,}\end{cases}\] and $n=dk+r$ is the Euclidean division of $n$ by $k.$
\end{defn}
In the same way, we define the type $B$ width-$k$ order polynomial of a poset $P$ having $p$ descents of type $B$, denoted by $\Omega_B(P,k,p)$, as the sum of the type $B$ width-$k$ order polynomials of its linear extensions
$$\Omega_B(P,k,p)= \displaystyle\sum_{\pi \in\mathfrak{L}(P) }\Omega_B(\pi,k,p),$$
where, for any permutation $\pi\in B_n$, the type $B$ width-$k$ order polynomial is
\begin{equation}\label{omega}
\Omega_B(\pi,k,p):=| \{f:[n-k+1] \to [p]/
\end{equation}
\begin{gather*}
1 \leq f(\pi(1))\leq f(\pi(1+k))\leq f(\pi(1+2 k))\leq \cdots f(\pi(1+\tau_d k))< \\ f(\pi(2))\leq f(\pi(2+k))\leq f(\pi(2+2 k))\leq \cdots f(\pi(2+\tau_d k))<\cdots < \\f(\pi(k))\leq f(\pi(k+k))\leq f(\pi(k+2k))\leq \cdots \leq f(\pi(k+\tau_d k)) \leq p,\\ and ~f(\pi(s))< f(\pi(s+k)),~ if~ s\in des^{B}_k(\pi)\}|.
\end{gather*}
In the following theorem, we give the analogous generating function
type $B$ polynomials order of Proposition 3.10 in \cite{26} in terms
of statistics of width-$k$ descents.
\begin{thm}\label{korder} For a given permutation $\pi \in B_n$, the generating function for width-$k$ type $B$ order polynomials
is given by
\begin{equation}
\displaystyle \sum_{p \geq 0} \Omega_B(\pi,k,p) x^{p}= \frac {x^{des^{B}_k(\pi)}}{(1-x)^{n-k+2}}.
\end{equation}
\end{thm}

\begin{dem}
Let $\mathfrak{L}(P)=\{\pi\}$, where $\pi$ has width-$k$ descents
counting an extra width-$k$ descent at the end (we assume
$\pi(k+\tau_d k)>\pi(n+1)$). Then, $\Omega_B(\pi;k;p)$ is the number
of solutions of Equation (\ref{omega}):
\begin{gather*}
1 \leq f(\pi(1))\leq f(\pi(1+k))\leq f(\pi(1+2 k))\leq \cdots f(\pi(1+\tau_d k))< \\ f(\pi(2))\leq f(\pi(2+k))\leq f(\pi(2+2 k))\leq \cdots f(\pi(2+\tau_d k))<\cdots < \\f(\pi(k))\leq f(\pi(k+k))\leq f(\pi(k+2k))\leq \cdots \leq f(\pi(k+\tau_d k)) \leq p-(des^{B}_k(\pi)-1),
\end{gather*}
which is equal to the number of ways choosing $n-k+1$ things from
$p-des^{B}_k(\pi)+1$ with repetitions. Therefore, the number is
$\binom{p-des^{B}_k(\pi)+1+n-k}{n-k+1}=\binom{p-des^{B}_k(\pi)+n-k+1}{n-k+1}.$
So, we obtain $$\displaystyle \sum_{p \geq 0} \Omega_B(\pi,k,p)
x^{p}=\displaystyle \sum_{p \geq 0}
\binom{p-des^{B}_k(\pi)+n-k+1}{n-k+1}x^{p}=\frac{x^{des^{B}_k(\pi)}}{(1-x)^{n-k+2}}~.$$
\end{dem}

Petersen \cite{10} gave a relation between enriched $P$-partitions
and quasisymmetric functions. Hence, we can see this relation in
terms of width-$k$ statistic which the connection helps us to prove
Theorem \ref{gammaB}.
The basic theory of enriched $P$-partitions is due to Stembridge \cite{22}. An enriched $P$-partitions and
 a left enriched $P$-partitions of type $A$ are defined as follow:
let $\mathbb{P'}$ denote the set of nonzero integers, totally
ordered given by$$-1<+1<-2<+2<-3<+3<\cdots,$$ and
$\mathbb{P^{(\ell)}}$ to be the integers with the following total
order
$$0<-1<+1<-2<+2<-3<+3<\cdots.$$
For any totally ordered
set $X=\{x_1, x_2, \ldots\}$, let $\mathbb{X'}$ and
$\mathbb{X^{(\ell)}}$ to be the sets given, respectively, by  $$\{-x_1, x_1, -x_2, \cdots
\}\quad and\quad \{x_0, -x_1, x_1, -x_2, \cdots \},$$ with total order $$
x_0 <-x_1 <x_1 <-x_2 <x_2 < \cdots .$$
\begin{defn} $[\cite{10},Definition 4.1]$ An enriched $P$-partition $($resp. left enriched $P$-partition$)$ is an order-preserving map $f:P \to \mathbb{X'}$ $($resp. $\mathbb{X^{(\ell)})}$ such that, for all $i<_P j$ in $P$, we have
\begin{enumerate}
\item if $i<_P j$ in $\mathbb{Z}$ then $f(i)\leq^{+} f(j)$,
\item if $i<_P j$ and $i>j$ in  $\mathbb{Z}$ then $f(i)\leq^{-} f(j).$
\end{enumerate}
\end{defn}
\noindent For type $B$ poset, an enriched $P$-partition of type $B$ differs from type $A$ only in the addition of a third condition, $f$($-i$) $= -f$($i$) (for more details, see [\cite{10}, Definition 4.8]).

\noindent Let $\varepsilon(P)$ be the set of all enriched $P$-partitions and $\varepsilon^{(\ell)}(P)$ be the set of left enriched $P$-partitions. The number of (left) enriched $P$-partitions is finite if we assume that $|X|=p$. In this case, we define
\emph{the enriched order polynomial}, denoted by $\Omega'(P,p)$, to be the number of enriched $P$-partitions $f:P\to X'$ and \emph{ the left enriched order polynomial}, denoted by $\Omega^{(\ell)}(P,p)$, to be the number of left enriched $P$-partitions $f:P\to X^{(\ell)}$.

\noindent Following Gessel \cite{21}, a quasisymmetric function is a map for which the coefficient of $x_{i_{1}}^{\alpha_1} x_{i_{2}}^{\alpha_2}\cdots x_{i_{p}}^{\alpha_p}$ is the same for all fixed tuples of integers $(\alpha_1, \alpha_2,\ldots,\alpha_p)$ and for all $i_{1}< i_{2}<\cdots<i_{p} $.
For any subset $D=\{d_1< d_2<\cdots <d_{p-1}\}$ of $[n]$, the quasisymmetric functions are characterized by two common bases: the monomial quasisymmetric functions, denoted $M_D$, and the fundamental quasisymmetric functions, denoted $F_D$, given by
\begin{eqnarray*}
M_D&=&\displaystyle \sum_{i_{1}<i_{2}<\cdots<i_{p}} x_{i_{1}}^{d_1} x_{i_{2}}^{d_2-d_1}\cdots x_{i_{p}}^{n-d_{p-1}} = \displaystyle \sum_{i_{1}<i_{2}<\cdots<i_{p}} x_{i_{1}}^{\alpha_1}  x_{i_{2}}^{\alpha_2} \cdots  x_{i_{p}}^{\alpha_p},\\
F_D &=&\displaystyle \sum_{D\subset T\subset[n-1]} M_T
=\displaystyle \sum_{{i_{1}\leq i_{2}\leq \cdots\leq i_{p} \atop d\in D \Rightarrow i_d < i_{d+1}}} \prod_{d=1}^{n} x_{i_{d}}
=\Gamma_A(\pi),
\end{eqnarray*}
where $\Gamma_A(\pi)$ is the generating function for the type $A$ $P$-partitions of a permutation with descent set $D$. It is possible to recover the order polynomial of $\pi$ by specializing:
\begin{equation}\label{quasiA}
\Omega_A(\pi,p)=\Gamma_A(\pi)(1^p).
\end{equation}
For all $D\subset[n-1]$, the functions $M_D$ and $F_D$ span the quasisymmetric functions of degree n, denoted $Qsym_n$. The ring of quasisymmetric functions is defined by $$ Qsym:=\displaystyle\bigoplus_{n\geq0} Qsym_n.$$
Moreover, the generating function for enriched $P$-partitions $f:P \to \mathbb{P'}$ is defined by $$\displaystyle\Delta_A(P):=\displaystyle\sum_{f\in\varepsilon(P)} \displaystyle\prod_{i=1}^{n} x_{|f(i)|}.$$
Evidently, $\displaystyle\Delta_A(P)$ is a quasisymmetric function and it satisfies
 $$\Omega'(P,p)=\displaystyle\Delta_A(P)(1^p).$$
Chow \cite{23} gave a connection between ordinary type $B$ $P$-partitions and type $B$ quasisymmetric functions. Furthermore, Petersen \cite{10} related the type $B$ quasisymmetric functions to left enriched $P$-partitions and type $B$ enriched $P$-partitions.\\
For a fixed $n$ and for any subset $D=\{d_1< d_2<\cdots<d_{p-1}\}$ of $\{0,1,\ldots, n\}$, the monomial $N_D$ and the fundamental quasisymmetric functions of type $B$, $L_D$, are defined, respectively, by
\begin{eqnarray*}
N_D &=&\displaystyle \sum_{0<i_{2}<\cdots<i_{p}} x_0^{d_1} x_{i_{2}}^{d_2-d_1}\cdots x_{i_{p}}^{n-d_{p-1}} = \displaystyle \sum_{0<i_{2}<\cdots<i_{p}} x_0^{\alpha_1}  x_{i_{2}}^{\alpha_2} \cdots  x_{i_{p}}^{\alpha_p},\\
L_D &=&\displaystyle \sum_{D\subset T\subset[0,n-1]} N_T
=\displaystyle \sum_{{0\leq i_{2}\leq \cdots\leq i_{p} \atop d\in D \Rightarrow i_d < i_{d+1}}} \prod_{d=1}^{n} x_{i_{d}}
=\Gamma_B(\pi),
\end{eqnarray*}
where $\Gamma_B(\pi)$ is the generating function for the ordinary type $B$ $P$-partitions, of any signed permutation with descent set $D$. In particular, we have
\begin{equation}\label{quasiB}
\Omega_B(\pi,p)=\Gamma_B(\pi)(1^{p+1}).
\end{equation}
Similar to type $A$, these functions form a basis for the type $B$ quasisymmetric functions of degree $n$, defined by
$$BQsym:=\displaystyle\bigoplus_{n\geq 0} BQsym_n.$$
The generating function for left enriched $P$-partitions $f$: $P \to \mathbb{P^{(\ell)}}$, is defined by $$\Delta^{(\ell)}(P):=\displaystyle\sum_{f\in \varepsilon^{(\ell)}(P)}\displaystyle\prod_{i=1}^{n} x_{|f(i)|}. $$
It is also clear that $\displaystyle\Delta^{(\ell)}(P)$ is a quasisymmetric function and we have $$\Omega^{(\ell)}(P,p)=\Delta^{(\ell)}(P)(1^{p+1}).$$
For any two subsets of integers $D$ ant $T$, define the set $D+k=\{d+k; d\in D\}$ with, $1\leq k \leq n$ and define the symmetric set difference by $$D\Delta T=(D\cup T)\backslash (D\cap T).$$
Using the standardization map $\varphi$ on $\symm_n$ defined in Section $2$, we can write the width-$k$ descent as in the following form $$des_k^B(\pi)=des^B(\pi_k)+\displaystyle\sum_{i=1}^{k-1} (des^B(\pi_i)-\epsilon(\pi_i)).$$
So, we count $0$ as a descent just in $\pi_k$, and all other permutations, we see it as a descent of type $A$. In addition, we can define $\Omega_B(D,k,p)$ by the ordinary type $B$ order polynomial of any signed permutation with width-$k$ descent set $D$ by
\begin{equation}\label{omegaB}
\Omega_B(D,k,p)=\Omega_B(D_k,p_k)\displaystyle \prod_{i=1}^{k-1} \Omega_B^{'}(D_i,p_i),
\end{equation}
with $\Omega_B(D_k,p_k)$ is the ordinary type $B$ order polynomial of signed permutation with descent set $D_k$ and $\Omega_B^{'}(D_i,p_i)$ is the ordinary type $B$ order polynomial of signed permutation with descent set $D_i\backslash\{0\}$,
  where $D_i$ is the set of descents in each $\pi_i$ and $p_i= p-k+i$, for any $1\leq i\leq k$. Thus,
   $$\displaystyle\bigcup_{i=1}^{k} D_i=D.$$

To study the left enriched $P$-partitions, $\Omega^{(\ell)}$($P,p$), it is sufficient to consider the case where $P$ is a permutation. So, it is possible to characterize the set of all left enriched $\pi$-partition in terms of descent set. Left peaks in $\symm_n$ are a special case of type $B$ peaks, which are naturally related to type $B$ descents. Thus, in terms of width-$k$ descent sets, we can define the set of all width-$k$ left enriched $\pi$-permutation, denoted by $\Omega^{(\ell)}$($\pi,k,p$), in which we consider $0$ as a descent only in $\pi_k$. Then, for $\pi \in \symm_n$, we can write
\begin{equation}\label{omegaleft}
\Omega^{(\ell)}(\pi,k,p)=\Omega^{(\ell)}(\pi_k,p_k)\prod_{i=1}^{k-1} \Omega^{(\ell)'}(\pi_i,p_i),
\end{equation}
where $\Omega^{(\ell)'}(\pi_i,p_i)$ is the left enriched order polynomial with $Des(\pi_i)\subset[\tau_d]$.
It is important to observe that $$Lpeak_k(\pi)=Lpeak(\pi_k)\displaystyle\bigcup_{i=1}^{k-1} Peak(\pi_i).$$ Thus, $$lpeak_k(\pi)=lpeak(\pi_k)+\displaystyle\sum_{i=1}^{k-1} peak(\pi_i).$$
The generating function for enriched $\pi$-partitions depends on the set of peaks. Moreover, the generating function for left enriched $\pi$-partitions depends on the set of left peaks. We can write $\Delta^{(\ell)}(\pi)$ according to the monomial and fundamental quasisymmetric functions of type $B$. In this case, by using the map $\varphi$ introduced in Section $2$, we can define the width-$k$ left enriched $P$-partitions by $$\Delta^{(\ell)}(\pi,k)=\Delta^{(\ell)}(\pi_k) \displaystyle\prod_{i=1}^{k-1} \Delta_A(\pi_i),$$
with nonnegative coefficients.
This coefficient is equal to the product number of enriched and left enriched $\pi$-partitions $f$ satisfying
$$(|f(\pi_1(1))|,|f(\pi_1(2))|,\ldots,|f(\pi_1(1+\tau_d k))|)=(1,\ldots1,\ldots,p_{1},\ldots p_{1}),$$ $$(|f(\pi_2(1))|,|f(\pi_2(2))|,\ldots,|f(\pi_2(1+\tau_d k))|)=(1,\ldots1,\ldots,p_{2},\ldots p_{2}),$$
$$\vdots$$
$$(|f(\pi_k(1))|,|f(\pi_k(2))|,\ldots,|f(\pi_k(1+d k))|)=(0,\ldots,0,1,\ldots1,\ldots,p_{k},\ldots p_{k}).$$
Applying Proposition $3.5$-\cite{22}, Theorem $6.6$-\cite{10} and referring to Equation (\ref{omegaB}), we obtain
$$\Delta^{(\ell)}(\pi,k)=2^{lpeak(\pi_k)} \displaystyle\sum_{ D_{k}\subset[0,d-1] \atop Lpeak(\pi_k)\subset D_{k}\Delta(D_{k}+1)} L_{D_k} \displaystyle \prod_{i=1}^{k-1} 2^{peak(\pi_i)+1} \displaystyle\sum_{ D_{i}\subset[\tau_d] \atop Peak(\pi_i)\subset D_{i}\Delta(D_{i}+1) }F_{D_i}.$$
We can also write
\begin{align*}
\Omega^{(\ell)}(\pi,k,p)&=2^{lpeak(\pi_k)}  \displaystyle\sum_{ D_{k}\subset[0,d-1] \atop Lpeak(\pi_k)\subset D_{k}\Delta(D_{k}+1)} \Omega_B(D_k,p_k) \displaystyle \prod_{i=1}^{k-1} 2^{peak(\pi_i)+1} \displaystyle\sum_{ D_{i}\subset[\tau_d] \atop Peak(\pi_i)\subset D_{i}\Delta(D_{i}+1)} \Omega'_B(D_i,p_i),\\
&=2^{lpeak_k(\pi)+k-1} \displaystyle\sum_{ D\subset[0,n-k] \atop Lpeak_k(\pi)\subset D\Delta(D+k)} \Omega_B(D,k,p).
\end{align*}
%From this result, we can find the following theorem.
\begin{thm}\label{thlpeak} The following generating function for width-$k$ left enriched order polynomials holds true
$$\displaystyle \sum_{p\geq 0} \Omega^{\mathcal{(\ell)}}(\pi,k,p) x^{p}=2^{k-1} \frac{(1+x)^{n-k+1}}{(1-x)^{n-k+2}} \left(\frac{4x}{(1+x)^{2}}\right)^{lpeak_k(\pi)}.$$
\end{thm}
\begin{dem} For any permutation $\pi$ in $\symm_n$, we have
$$\displaystyle \sum_{p\geq 0} \Omega^{(\ell)}(\pi,k,p)x^p=\displaystyle \sum_{p\geq 0}2^{lpeak_k(\pi)+k-1} \displaystyle\sum_{ D\subset[0,n-k] \atop Lpeak_k(\pi)\subset D\Delta(D+k)} \Omega_B(D,k,p)x^p.$$
Applying the generating function for width-k type $B$ order polynomials defined in Theorem \ref{korder}, we obtain
\begin{align*}
\displaystyle \sum_{p\geq 0} \Omega^{(\ell)}(\pi,k,p)x^p = \frac{2^{lpeak_k(\pi)+k-1}}{(1-x)^{n-k+2}} \displaystyle\sum_{ D\subset[0,n-k] \atop Lpeak_k(\pi)\subset D\Delta(D+k)} x^{|D|}.
\end{align*}
It is not complicated to specify the generation function for $D$ sets by size. For all $s\in Lpeak_k$($\pi$), we have precisely $s$ or $s-k$ is in $D$. So, there is still $n-k+1-2lpeak_k$($\pi$) elements in $\{0,1,\ldots, n-k\}$ that can be contained in $D$ or not. Thus, we obtain
\begin{align*}
\displaystyle \sum_{ D\subset[0,n-k] \atop Lpeak_k(\pi)\subset D\Delta(D+k)} x^{|D|}&=\underbrace{(x+x)(x+x)\cdots(x+x)}_{lpeak_k(\pi)} \underbrace{(1+x)(1+x)\cdots(1+x)}_{n-k+1-2lpeak_k(\pi)}\\
&=(2x)^{lpeak_k(\pi)} (1+x)^{n-k+1-2lpeak_k(\pi)}.
\end{align*}
By combining them all together, we deduce the desired identity
\begin{align*}
\displaystyle \sum_{p\geq 0} \Omega^{\mathcal{(\ell)}}(\pi,k,p) x^{p} &=\frac{2^{lpeak_k(\pi)+k-1}}{(1-x)^{n-k+2}}
(2x)^{lpeak_k(\pi)}(1+x)^{n-k+1-2lpeak_k(\pi)}\\
&=2^{k-1} \frac{(1+x)^{n-k+1}}{(1-x)^{n-k+2}} \left(\frac{4x}{(1+x)^{2}}\right)^{lpeak_k(\pi)}.
\end{align*}
\end{dem}

\noindent Recall that the number of permutations of $\symm_n$ with $p$ width-$k$ left peaks is $|\Gamma^{\mathcal{(\ell)}}_{  n,k,p}|$. So, the width-$k$ left peak polynomial can be defined as follows
\begin{equation}\label{left}
W_{n,k}^{(\ell)}(x)=\displaystyle \sum_{\pi \in \symm_n}x^{lpeak_k(\pi)} =\displaystyle \sum_{p=0}^{\lfloor \frac{n-k+1}{2}\rfloor} |\Gamma^{\mathcal{(\ell)}}_{  n,k,p}| x^p.
\end{equation}
The following identity for the type $B$ Eulerian polynomials is due to Stembridge, Proposition $7.1.b$-\cite{24},
  $$\displaystyle \sum_{p\geq0}(2p+1)^n x^p=\frac{\mathfrak{B}_n(x)}{(1-x)^{n+1}}.$$
In the following result, we establish the relation between the width-$k$ left peak polynomials and the width-$k$ Eulerian polynomials
of type $B$.
\begin{prop} The following identity holds true:
\begin{equation}
W_{n,k}^{(\ell)}\left(\frac{4x}{(1+x)^{2}}\right)=\frac{\mathfrak{WB}_{n,k}(x)}{2^{k-1}(1+x)^{n-k+1}}.
\end{equation}
\end{prop}
\begin{dem} Using Equation (\ref{omegaleft}), the number $\mathcal{N}$ of width-$k$ left enriched $P$-partitions $f$ is given by $$\mathcal{N}=(2p+1)^d~\left((2p+1)^{(d-1)}\right)^{r}~\left((2p+1)^{d}\right)^{k-r-1}.$$
Remark that $\pi_k\in \symm_d$, $(\pi_1,\ldots,\pi_{k-1})\in \symm_{d+1}^{r}\times \symm_{d}^{k-r-1}$ and the number of $f:[n-k+1]\to [p]^{(\ell)}$ is $(2p+1)^{n-k+1}$. Consequently, $\Omega^{(\ell)}(\pi,k,p)=(2p+1)^{n-k+1}$.
Let $P_B$ be a type $B$ poset whose elements are $\{-n,-n+1,\ldots,-1,1,\ldots,n\}$. In the present case, $P_B$ is an antichain of $[-n+k-1,n-k+1]$. Thus, the order polynomial $\Omega_B(\pi,k,p)$ is the same as of $\Omega^{(\ell)}(\pi,k,p)$.
On the other hand, we have $$\displaystyle \sum_{p\geq0}(2p+1)^{n-k+1} x^p=\frac{\mathfrak{WB}_{n,k}(x)}{(1-x)^{n-k+2}}.$$
Using Theorem \ref{thlpeak}, we obtain
\begin{align*}
2^{k-1} \frac{(1+x)^{n-k+1}}{(1-x)^{n-k+2}} W_{n,k}^{(\ell)}\left(\frac{4x}{(1+x)^{2}}\right)&=\displaystyle \sum_{p\geq0} (2p+1)^{n-k+1} x^p\\
&=\frac{\mathfrak{WB}_{n,k}(x)}{(1-x)^{n-k+2}}.
\end{align*}
The reorganization of the terms gives the desired identity, that is
\begin{equation*}\label{leftp}
2^{k-1}W_{n,k}^{(\ell)}\left(\frac{4x}{(1+x)^{2}}\right)=\frac{\mathfrak{WB}_{n,k}(x)}{(1+x)^{n-k+1}}.
\end{equation*}
\end{dem}

Now, we are able to prove Theorem \ref{gammaB}.\\

\noindent $Proof \; of\; Theorem\; \ref{gammaB}.$ Substituting $x$ by $\frac{4x}{(1+x)^2}$ in Equation (\ref{left}) and referring to Equation (\ref{leftp}), we obtain
$$\displaystyle\sum_{p=0}^{\lfloor \frac{n-k+1}{2}\rfloor} \gamma^{B}_{ n,k,p}~ x^{p} (1+x)^{n-k+1-2 p}=\displaystyle \sum_{p=0}^{\lfloor \frac{n-k+1}{2}\rfloor} 2^{2p+k-1} |\Gamma^{\mathcal{(\ell)}}_{n,k,p}| x^{p} (1+x)^{n-k+1-2 p}$$
which allows to conclude. \hfill$\blacksquare$
\\

\noindent For $1 \leq n \leq 6,~1\leq k \leq n$ and $0 \leq p \leq \lfloor \frac{n-k+1}{2}\rfloor$, we record a few values of $\gamma^{B}_{ n,k,p}$ in the following table.\newpage
\begin{table}[!h]
\centering
\begin{tabular}{|c|c|c c cc|}
\hline
\scriptsize n &\scriptsize k &
\multicolumn{4}{c|}{\centering{\scriptsize p}}\\[-.1cm]
\cline{3-6}
& &\scriptsize0&\scriptsize1&\scriptsize2&\scriptsize3\\[-.1cm]
\hline
\scriptsize1&\scriptsize1&\scriptsize1& & &\\[-.1cm]
\hline
\scriptsize2&\scriptsize1&\scriptsize1&\scriptsize4 & &\\[-.2cm]
 &\scriptsize2&\scriptsize4& & &\\[-.1cm]
\hline
\scriptsize3&\scriptsize1&\scriptsize1&\scriptsize20 & &\\[-.2cm]
 &\scriptsize2&\scriptsize12&\scriptsize0 & &\\[-.2cm]
 &\scriptsize3&\scriptsize24& & &\\[-.1cm]
\hline
\scriptsize4&\scriptsize1&\scriptsize1&\scriptsize 72&\scriptsize 0&\\[-.2cm]
 &\scriptsize2&\scriptsize24 &\scriptsize96 &&\\[-.2cm]
 &\scriptsize3&\scriptsize96&\scriptsize0 & &\\[-.2cm]
 &\scriptsize4&\scriptsize192& & &\\[-.1cm]
\hline
\scriptsize5&\scriptsize1&\scriptsize1& \scriptsize232&\scriptsize976 &\\[-.2cm]
 &\scriptsize2&\scriptsize80& \scriptsize480&\scriptsize 640&\\[-.2cm]
 &\scriptsize3&\scriptsize480&\scriptsize 0& &\\[-.2cm]
 &\scriptsize4&\scriptsize960&\scriptsize0 & &\\[-.2cm]
 &\scriptsize5&\scriptsize1920& & &\\[-.1cm]
\hline
\scriptsize6&\scriptsize1&\scriptsize1&\scriptsize716 &\scriptsize7664 &\scriptsize3904\\[-.2cm]
 &\scriptsize2&\scriptsize160&\scriptsize3520 &\scriptsize6400 &\\[-.2cm]
 &\scriptsize3&\scriptsize1440&\scriptsize 5760&\scriptsize 0&\\[-.2cm]
 &\scriptsize4&\scriptsize5760&\scriptsize0 & &\\[-.2cm]
 &\scriptsize5&\scriptsize11520&\scriptsize 0& &\\[-.2cm]
 &\scriptsize6&\scriptsize23040& & &\\
\hline
\end{tabular}
\caption{The first few values of $\gamma^{B}_{  n,k,p}$.}
\end{table}
\section{Width-$k$ Eulerian polynomials of type $D$}
 In this section, we aim to study the $\gamma$-positivity of the width-$k$ Eulerian polynomials on the set $D_n$.\\

 Borowiec and  Mlotkowski \cite{12} have introduced a new array of type $D$ Eulerian numbers. They found, in particular, a recurrence relation for this array. Our purpose is to generalize these numbers and to give a new extensions of the Eulerian polynomials of type $D$, $D_n$($x$), and its $\gamma$-positivity with the statistics \emph{width-$k$ descent}. The case $k=1$ corresponds to the Eulerian polynomials $D_n$($x$) in the sense of Borowiec and  Mlotkowski \cite{12}.\\

 \noindent We start by defining the width-$k$ Eulerian polynomials of type $D$ as follows
 $$\mathfrak{WD}_{n,k}(x):=  F^{des^{D}_k}_n(x)= \displaystyle \sum_{\pi \in D_n} x^{des^{D}_k(\pi)}.$$
The following theorem is the mean result of this section.
\begin{thm}Let $n\geq 1$. Then, the following statements holds true:
\begin{enumerate}
\item If $k=1$, then we have
$$\mathfrak{WD}_{n,1}(x) = D_n(x),$$
\item If $2 \leq k \leq n $, then we have
\begin{equation}\label{desd}
\mathfrak{WD}_{n,k}(x) = 2^{n-d-1} M_{n,k} B_{d}(x) A_{d}^{k-r-1}(x) A_{d+1}^{r}(x).
\end{equation}
\end{enumerate}
\end{thm}
\begin{dem} If $2\leq k \leq n$, we use the standardization map $\psi$ defined in Section $2$ and restrained to the set $D_n$. So, the desired identity follows completely in the same way as the identity (\ref{desb}) of Theorem \ref{thmb}, by applying the same reasoning and taking into account that $|D_n|=\frac{|B_n|}{2}$.

\end{dem}

\vspace{.2cm}
Denote $\bar{D_n}=B_n \backslash D_n$ and,
$$\mathrm{WD}_{n,k,p}=\{\pi \in D_n; des^{D}_k(\pi)=p\},$$
$$\mathrm{W\bar{D}}_{n,k,p}=\{\pi \in \bar{D}_{n}; des^{\bar{D}}_k(\pi)=p\}.$$
So, $\mathrm{WD}_{n,k,p}=\mathrm{WB}_{n,k,p}\bigcap D_n$ and $\mathrm{W\bar{D}}_{n,k,p}=\mathrm{WB}_{n,k,p} \backslash  D_n$. The cardinals of these sets will be denoted by $d(n,k,p)$ and $\bar{d}(n,k,p)$, respectively. Since $\mathrm{WB}_{n,k,p}=\mathrm{WD}_{n,k,p} \bigcup \mathrm{W\bar{D}}_{n,k,p} $, we have
\begin{equation}\label{d+d} b(n,k,p)=d(n,k,p)+ \bar{d}(n,k,p).\end{equation}

\noindent Now, we give the following result in which we generalize Proposition $4.1$-\cite{12} (for $k=1)$.
\begin{prop} Let $0 \leq p \leq n$. Then, the following statements
are satisfied:
 \begin{enumerate}
\item If $n$ is even with $1 \leq k \leq n$ or $n$ is odd with $2 \leq k
\leq n $, then
\begin{equation}\label{dd}\delta(n,k,p)=\delta(n,k,n-k+1-p)~ and
~\bar{\delta}(n,k,p)=\bar{\delta}(n,k,n-k+1-p).
\end{equation}
\item If $n$ is odd with $k=1$, then
$$\delta(n,1,p)=\bar{\delta}(n,1,n-p)~ and ~ \bar{\delta}(n,1,p)=\delta(n,1,n-p).$$
\end{enumerate}
\end{prop}
\begin{dem} Let $\xi: B_n  \to  B_n $ be a map satisfying $\xi(\pi)=(-\pi(1), \cdots, -\pi(n))$, for any $\pi=(\pi(1), \cdots, \pi(n)) \in B_n$.
 It is easy to see that $\xi$ is bijective from $\mathrm{WD}_{n,k,p}$ onto ${WD}_{n,k,n-k+1-p}$ and from $\mathrm{W\bar{D}}_{n,k,p}$
 onto $\mathrm{W\bar{D}}_{n,k,n-k+1-p} $. So, Equations (\ref{dd}) hold true if $n$ is even with $1 \leq k \leq n$. Moreover,
 $\xi$ is bijective from
 $\mathrm{WD}_{n,k,p}$ onto $\mathrm{W\bar{D}}_{n,k,n-k+1-p}$ if $n$ is odd with $2 \leq k \leq n$. By using Eqation (\ref{d+d}), we obtaine the desired claim.

\end{dem}

\noindent For $k=1$, Borowiec and Mlotkowski \cite{12} showed the
following recurrence relations.
\begin{prop}
Let $0\leq p \leq n$. Then, the following identities hold true:
\begin{eqnarray*}
d(n,1,p)&=&(2p+1) d(n-1,1,p)+(2n-2p+1) d(n-1,1,p-1)+
(-1)^{p}\binom{n-1}{p-1},\\
\bar{d}(n,1,p)&=&(2p+1) \bar{d}(n-1,1,p)+(2n-2p+1)
\bar{d}(n-1,1,p-1)- (-1)^{p}\binom{n-1}{p-1}. \end{eqnarray*}
\end{prop}
Referring to Equation (\ref{desd}), we give in the table \ref{delta}
some coefficients of the product polynomials $ B_{d}(x)
A_{d}^{k-r-1}(x) A_{d+1}^{r}(x)$, denoted by $\delta$($n,k,p$), for
$1\leq n \leq 6,~1< k \leq n$ and $0 \leq p \leq n-k+1$. Moreover,
in the table \ref{bardelta}, we give the first values of
$\bar{\delta}(n,k,p)$ .

\begin{table}[!h]
\begin{center}
\begin{tabular}{|c|c|c c c c c c c|}
\hline
\scriptsize n & \scriptsize k &
\multicolumn{7}{c|}{\centering{\scriptsize p}}\\[-.1cm]
\cline{3-9}
& & \scriptsize0 &\scriptsize1 &\scriptsize2 & \scriptsize3 &\scriptsize4 &\scriptsize 5 &\scriptsize6\\[-.1cm]
\hline
\scriptsize1&\scriptsize1&\scriptsize1&\scriptsize1 &  & & & & \\[-.1cm]
\hline
\scriptsize2&\scriptsize1&\scriptsize1&\scriptsize2 &\scriptsize1  & & & & \\[-.1cm]
 &\scriptsize2&\scriptsize1&\scriptsize1 &  & & & & \\[-.1cm]
\hline
\scriptsize3&\scriptsize1&\scriptsize1&\scriptsize10 & \scriptsize13 &\scriptsize 0& & & \\[-.1cm]
 &\scriptsize2&\scriptsize1&\scriptsize2 & \scriptsize1 & & & & \\[-.1cm]
 &\scriptsize3&\scriptsize1&\scriptsize1 &  & & & & \\[-.1cm]
\hline
\scriptsize4&\scriptsize1&\scriptsize1&\scriptsize36 &\scriptsize 118 & \scriptsize36&\scriptsize1 & & \\[-.1cm]
 &\scriptsize2&\scriptsize1&\scriptsize7 & \scriptsize7 &\scriptsize1 & & & \\[-.1cm]
 &\scriptsize3&\scriptsize1&\scriptsize2 & \scriptsize1 & & & & \\[-.1cm]
 &\scriptsize4&\scriptsize1&\scriptsize1 &  & & & & \\[-.1cm]
\hline
\scriptsize5&\scriptsize1&\scriptsize1&\scriptsize116 &\scriptsize846  &\scriptsize836 &\scriptsize 121& \scriptsize0& \\[-.1cm]
 &\scriptsize2&\scriptsize1&\scriptsize10 & \scriptsize26 & \scriptsize10&\scriptsize 1& & \\[-.1cm]
 &\scriptsize3&\scriptsize1&\scriptsize3 &\scriptsize 3 &\scriptsize 1& & & \\[-.1cm]
 &\scriptsize4&\scriptsize1&\scriptsize2 & \scriptsize 1& & & & \\[-.1cm]
 &\scriptsize5&\scriptsize1&\scriptsize1 &  & & & & \\[-.1cm]
\hline
\scriptsize6&\scriptsize1&\scriptsize1&\scriptsize358 &\scriptsize5279  &\scriptsize11764 &\scriptsize5279 &\scriptsize358 & \scriptsize1\\[-.1cm]
 &\scriptsize2&\scriptsize1&\scriptsize27 &\scriptsize116  & \scriptsize116& \scriptsize27& \scriptsize1& \\[-.1cm]
 &\scriptsize3&\scriptsize1&\scriptsize8 & \scriptsize14 &\scriptsize\scriptsize8 &\scriptsize 1& & \\[-.1cm]
 &\scriptsize4&\scriptsize1&\scriptsize3&\scriptsize 3 &\scriptsize1 & & & \\[-.1cm]
 &\scriptsize5&\scriptsize1&\scriptsize2 & \scriptsize1 & & & & \\[-.1cm]
 &\scriptsize6&\scriptsize1&\scriptsize1 &  & & & & \\
\hline
\end{tabular}
\caption{The first values of $\delta(n,k,p)$.}\label{delta}

\begin{tabular}{|c|c|c c c c c c c|}
\hline
\scriptsize n & \scriptsize k &
\multicolumn{7}{c|}{\centering{\scriptsize p}}\\[-.1cm]
\cline{3-9}
& & \scriptsize0 &\scriptsize1 &\scriptsize2 & \scriptsize3 &\scriptsize4 &\scriptsize 5 &\scriptsize6\\[-.1cm]
\hline
\scriptsize1&\scriptsize1&\scriptsize0&\scriptsize1 &  & & & & \\[-.1cm]
\hline
\scriptsize2&\scriptsize1&\scriptsize0&\scriptsize4 &\scriptsize0  & & & & \\[-.1cm]
 &\scriptsize2&\scriptsize1&\scriptsize1 &  & & & & \\[-.1cm]
\hline
\scriptsize3&\scriptsize1&\scriptsize0&\scriptsize13 & \scriptsize10 &\scriptsize 1& & & \\[-.1cm]
 &\scriptsize2&\scriptsize1&\scriptsize2 & \scriptsize1 & & & & \\[-.1cm]
 &\scriptsize3&\scriptsize1&\scriptsize1 &  & & & & \\[-.1cm]
\hline
\scriptsize4&\scriptsize1&\scriptsize0&\scriptsize40 & \scriptsize112 &\scriptsize 40&\scriptsize0 & & \\[-.1cm]
 &\scriptsize2&\scriptsize1&\scriptsize7 & \scriptsize7 &\scriptsize1 & & & \\[-.1cm]
 &\scriptsize3&\scriptsize1&\scriptsize2 & \scriptsize1 & & & & \\[-.1cm]
 &\scriptsize4&\scriptsize1&\scriptsize1 &  & & & & \\[-.1cm]
\hline
\scriptsize5&\scriptsize1&\scriptsize0&\scriptsize121 &\scriptsize836  &\scriptsize846 & \scriptsize116& \scriptsize1& \\[-.1cm]
 &\scriptsize2&\scriptsize1&\scriptsize10 &\scriptsize 26 & \scriptsize10&\scriptsize 1& & \\[-.1cm]
 &\scriptsize3&\scriptsize1&\scriptsize3 & \scriptsize3 & \scriptsize1& & & \\[-.1cm]
 &\scriptsize4&\scriptsize1&\scriptsize2 & \scriptsize 1& & & & \\[-.1cm]
 &\scriptsize5&\scriptsize1&\scriptsize1 &  & & & & \\[-.1cm]
\hline
\scriptsize6&\scriptsize1&\scriptsize0&\scriptsize364 &\scriptsize5264  &\scriptsize11784 &\scriptsize5264 &\scriptsize364 & \scriptsize0\\[-.1cm]
 &\scriptsize2&\scriptsize1&\scriptsize27 &\scriptsize116  &\scriptsize 116&\scriptsize 27&\scriptsize 1& \\[-.1cm]
 &\scriptsize3&\scriptsize1&\scriptsize8 &\scriptsize 14 &\scriptsize8 & \scriptsize1& & \\[-.1cm]
 &\scriptsize4&\scriptsize1&\scriptsize3&\scriptsize 3 &\scriptsize1 & & & \\[-.1cm]
 &\scriptsize5&\scriptsize1&\scriptsize2 & \scriptsize1 & & & & \\[-.1cm]
 &\scriptsize6&\scriptsize1&\scriptsize1 &  & & & & \\
\hline
\end{tabular}
\caption{The first values of $\bar{\delta}(n,k,p)$.}\label{bardelta}
\end{center}
\end{table}

Remark that, for $0 \leq p \leq n-k+1$, with $k$ being the smallest
positive integer such that $n+1=2k+r$ $(0\leq r <k)$, we obtain
 the following recurrence relations
%Let $k$ be the smallest divisor that gives $2$ as quotient in Euclidean division of $n+1$ by $k$, means $n+1=2k+r$ for some $n,r\in \mathbb{N}$ with $0\leq r <k.$ Then, the number $d(n,k,p)$ admit the following recurrence relation
\begin{eqnarray*}
\delta(n,k,p)&=&\delta(n-2,k-1,p-1)+\delta(n-2,k-1,p),\quad n\geq 3,\\
\bar{\delta}(n,k,p)&=&\bar{\delta}(n-2,k-1,p-1)+\bar{\delta}(n-2,k-1,p),\quad
n \geq 4, \end{eqnarray*}
%Let $n$ and $k$ be positive integers for which $n+1=dk+r$ for some $d,r~\in \mathbb{N}$ with $0\leq r <k.$ Then, for k is the smallest divisor that gives us 2 as quotient. The numbers $d(n,k,p)~and~\bar{d}(n,k,p)$ admit the following recurrence relations
%$$d(n,k,p)=d(n-2,k-1,p-1)+d(n-2,k-1,p),$$ $$\bar{d}(n,k,p)=\bar{d}(n-2,k-1,p-1)+\bar{d}(n-2,k-1,p),$$
where,
\begin{eqnarray*}
\delta(n, k, 0)&=&\bar{\delta}(n, k, 0)=\delta(n, k,
n-k+1)=\bar{\delta}(n, k, n-k+1)= 1,~\textnormal{for}~2\leq k\leq
n,\\
\delta(n,k,-1)&=&\bar{\delta}(n,k,-1)=0,\\
\delta(n,1,0)&=&1,~~\bar{\delta}(n,1,0)=0
\end{eqnarray*}
 and
\[\delta(n,1,n)=\begin{cases} 1, & \text{if n is even;}\\0, & \text{if n is odd,}\end{cases},\qquad\bar{\delta}(n,1,n)=\begin{cases} 0, & \text{if n is even;}\\1, & \text{if n is odd.}\end{cases}\]
\begin{prob}
Is it possible to find the recurrence relations of  $\delta(n,k,p)$ and $\overline{\delta}(n,k,p)$ for all $1 \leq k \leq n?$
\end{prob}
\begin{rem}Let $n$ be a positive integer. Then, the following
statements hold:
\begin{enumerate}
\item If $n$ is even, then we obtain
$$deg(\mathfrak{WD}_{n,k}(x))=n-k+1~~for~all~1\leq k \leq n.$$
\item If $n$ is odd, then we obtain
\[deg(\mathfrak{WD}_{n,k}(x)) =\begin{cases} n-1, & \text{if $k=1;$}\\n-k+1, & \text{if $2 \leq k \leq n.$}\end{cases}\]
\end{enumerate}
Moreover, according to Table \ref{delta}, we can deduce that if $n$
is odd, with $k = 1$, then the polynomial $\mathfrak{WD}_{n,k}(x)$
is not unimodal and not symmetric. For example, we have
$$\mathfrak{WD}_{5,1}(x)= 1+ 116 x+ 846 x^{2}+836 x^{3}+121 x^{4}.$$
\end{rem}
\begin{thm}\label{gammaD} Let $n,\,k$ be positive integers with $2\leq k\leq n,$ and $n$ is not odd.
Then, the following identity holds true
\[\mathfrak{WD}_{n,k}(x)=\displaystyle \sum_{\pi \in D_n} x^{des^{D}_k(\pi)}=\displaystyle\sum_{p=0}^{\lfloor \frac{n-k+1}{2}\rfloor}
 \frac{\gamma^{B}_{ n,k,p}}{2}~ x^{p} (1+x)^{n-k+1-2 p},\].
\end{thm}
\begin{dem}
If $n$ is odd and $k=1$, there are no permutations $\pi$ in $D_n$
such that $des^{D}_1(\pi)=n$. Thus, $\mathfrak{WD}_{n,1}(x)$ is not symmetric, and therefore $\mathfrak{WD}_{n,1}(x)$ is not $\gamma$-positive.\\
If $n$ is even, then for all $1 \leq k \leq n$, the number of permutations whose $des^{D}_k(\pi)=n-k+1-p$ is equal to the number of permutations of
 which $des^{D}_k(\pi)=p$ with $0 \leq p \leq \lfloor \frac{n-k+1}{2}\rfloor$. Thus, $\mathfrak{WD}_{n,k}(x)$ is symmetric with center of symmetry $\lfloor\frac{n-k+1}{2}\rfloor$. Moreover, it is easy to see that $\mathfrak{WD}_{n,k}(x)$ is unimodal.\\
 Therefore, since for all $2\leq k \leq n $ the width-$k$ Eulerian polynomials of type $D$ is $\gamma$-positive, and the cardinal of this group
  is equal to $\frac{|B_{n,k}|}{2}$, thus $\gamma^{D}_{ n,k,p}=\frac{\gamma^{B}_{ n,k,p}}{2}.$

\end{dem}

\noindent For $2 \leq n \leq 6,~1\leq k \leq n$ and $0 \leq p \leq
\lfloor \frac{n-k+1}{2}\rfloor$, we give in Table \ref{deltaD} some
values of $\gamma^{D}_{ n,k,p}$.
\begin{table}[!!h]
\centering
\begin{tabular}{|c|c|c c c c|}
\hline
\scriptsize n &\scriptsize k &
\multicolumn{4}{c|}{\centering{\scriptsize p}}\\[-.1cm]
\cline{3-6}
& &\scriptsize0&\scriptsize1&\scriptsize2&\scriptsize3\\[-.1cm]
\hline
\scriptsize2&\scriptsize1&\scriptsize1&\scriptsize0 & &\\[-.1cm]
 &\scriptsize2&\scriptsize2& & &\\[-.1cm]
\hline
\scriptsize3&\scriptsize1& & & &\\[-.1cm]
 &\scriptsize2&\scriptsize6&\scriptsize0 & &\\[-.1cm]
 &\scriptsize3&\scriptsize12& & &\\[-.1cm]
\hline
\scriptsize4&\scriptsize1&\scriptsize1& \scriptsize32&\scriptsize 48&\\[-.1cm]
 &\scriptsize2&\scriptsize12&\scriptsize48 & &\\[-.1cm]
 &\scriptsize3&\scriptsize48&\scriptsize0 & &\\[-.1cm]
 &\scriptsize4&\scriptsize96& & &\\[-.1cm]
\hline
\scriptsize5&\scriptsize1& &  &  &\\[-.1cm]
 &\scriptsize2&\scriptsize40& \scriptsize240&\scriptsize 320&\\[-.1cm]
 &\scriptsize3&\scriptsize240&\scriptsize 0& &\\[-.1cm]
 &\scriptsize4&\scriptsize480&\scriptsize0 & &\\[-.1cm]
 &\scriptsize5&\scriptsize960& & &\\[-.1cm]
\hline
\scriptsize6&\scriptsize1&\scriptsize1&\scriptsize352 &\scriptsize3856 &\scriptsize1920\\[-.1cm]
 &\scriptsize2&\scriptsize80&\scriptsize1760 &\scriptsize3200 &\\[-.1cm]
 &\scriptsize3&\scriptsize720& \scriptsize2880&\scriptsize 0&\\[-.1cm]
 &\scriptsize4&\scriptsize2880&\scriptsize0 & &\\[-.1cm]
 &\scriptsize5&\scriptsize5760& \scriptsize0& &\\[-.1cm]
 &\scriptsize6&\scriptsize11520& & &\\
\hline
\end{tabular}
\caption{The first values of $\gamma^{D}_{ n,k,p}.$}\label{deltaD}
\end{table}
\begin{prob}
Is it possible to find the recurrence relation of $\gamma^{D}_{ n,1,p}$, if $n$ is even ($k=1$)?
\end{prob}
\section*{Acknowledgments}
We would like to thank Prof. Jaing Zeng for reading an earlier
version of this manuscript and for their comments.

\end{document}